\documentclass[11pt,psamsfonts]{amsart}
 \usepackage{amsmath}
 \usepackage{amssymb}
 \usepackage{amscd}
\usepackage[all]{xy}
\usepackage[colorlinks=true, pdfstartview=FitV, linkcolor=blue, citecolor=blue, urlcolor=blue]{hyperref}

\newtheorem*{introtheorem1}{Slice Theorem}
\newtheorem*{introtheorem2}{Quotient Theorem}
\newtheorem*{introtheorem3}{Stratification Theorem}
\newtheorem*{theoremproper}{Theorem}
\newtheorem*{theoremcompact}{Theorem}

\newtheorem*{morseineq}{Morse inequalities}

\numberwithin{equation}{subsection}
\newtheorem{theorem}[subsection]{Theorem}
\newtheorem{lemma}[subsection]{Lemma}
\newtheorem{proposition}[subsection]{Proposition}
\newtheorem{corollary}[subsection]{Corollary}

\theoremstyle{definition}

\theoremstyle{remark}
\newtheorem{remark}[subsection]{Remark}

\newtheorem{example}[subsection]{Example}

\def\inv{^{-1}}

\def\CC{\mathbb C}
\def\PP{\mathbb P}

\def\quot#1#2{{#1/\!\!/#2}}
\def\twist#1#2#3{#1\times^{#2}#3}
\def\c{^\CC}

\def\lie#1{\mathfrak{ #1}}
\def\lieh{\lie h}
\def\liek{\lie k}
\def\liep{\lie p}
\def\lieq{\lie q}
\def\liea{\lie a}
\def\lier{\lie r}
\def\lieg{\lie g}
\def\liem{\lie m}
\def\lieu{\lie u}

\def\grad{\operatorname{grad}}
\def\RR{{\mathbb R}}

\newcommand{\Ad}{\operatorname{Ad}}
\newcommand{\ad}{\operatorname{ad}}
\def\muip{\mu_{\im\liep}}
\def\mup{\mu_{\liep}}

\def\mua{\mu_{\liea}}

\def\mip{\mathcal M_{\im\liep}}
\def\mp{\mathcal M_{\liep}}

\def\cp{\mathcal C_{\liep}}

\def\bp{\mathcal B_{\liep}}

\def\m{\mathcal M}
\def\im{\mathrm{i}}
\def\sgmip{\mathcal{S}_G(\mip)}

\def\SS{\mathcal{S}}

\def\Id{\operatorname{Id}}

\def\norm#1{\lVert#1\rVert}
\def\im{\mathrm{i}}

\newcommand{\co}[1]{H^n(#1)}
\newcommand{\eco}[2]{H^n_{#1}(#2)}
\newcommand{\poin}[2]{P_t^{#1}(#2)}
\newcommand{\ce}{\mathbb{C}}
\newcommand{\zet}{\mathbb{Z}}
\newcommand{\qu}{\mathbb{Q}}
\newcommand{\er}{\mathbb{R}}
\newcommand{\N}{\mathbb{N}}
\newcommand{\peinsce}{\mathbb{P}_1(\ce)}
\newcommand{\pzweice}{\mathbb{P}_2(\ce)}
\newcommand{\pzweier}{\mathbb{P}_2(\er)}
\newcommand{\Sldreice}{\operatorname{SL}_3(\ce)}
\newcommand{\Sldreier}{\operatorname{SL}_3(\er)}
\newcommand{\slzweier}{\operatorname{SL}_2(\er)}
\newcommand{\slzweice}{\operatorname{SL}_2(\ce)}
\newcommand{\sozweier}{\operatorname{SO}_2(\er)}
\newcommand{\acts}{\mbox{\raisebox{0.26ex}{\;\tiny{$\bullet$}\;}}}

\title[Stratifications]{Stratifications with respect to actions of real reductive groups}
\author{Peter Heinzner}
\email{heinzner@cplx.rub.de}
\address{Fakult\"at f\"ur Mathematik\\
Ruhr-Universit\"at Bochum\\
Universit\"atsstrasse 150\\
D - 44780 Bochum}
\author{Gerald W. Schwarz}
\email{schwarz@brandeis.edu}
\address{Department of Mathematics\\
Brandeis University\\
PO Box 549110\\
Waltham, MA 02454-9110}
\author{Henrik St\"otzel}
\email{henrik.stoetzel@rub.de}
\address{Fakult\"at f\"ur Mathematik\\
Ruhr-Universit\"at Bochum\\
Universit\"atsstrasse 150\\
D - 44780 Bochum}
\thanks{First author is partially supported by the Sonderforschungsbereich
SFB/TR12 of the Deutsche Forschungsgemeinschaft.\\
Second author is partially supported by NSA Grant
H98230-06-1-0023.\\
Third author is supported by the Sonderforschungsbereich SFB/TR12
of the Deutsche Forschungsgemeinschaft.}



\begin{document}

\begin{abstract}
We study the action of a real reductive group $G$ on a real
submanifold $X$ of a K\"ahler manifold $Z$. We suppose that the
action of $G$ extends holomorphically to an  action of the
complexified group $G\c$ and that with respect to a compatible
maximal compact subgroup $U$ of $G\c$ the action on $Z$ is
Hamiltonian. There is a corresponding gradient map $\mup\colon
X\to\lie p^*$ where $\lie g=\lie k\oplus\lie p$ is a Cartan
decomposition of $\lieg$. We obtain a Morse like function
$\eta_{\liep}:=\lVert\mup\rVert^2$ on $X$. Associated
to critical points of $\eta_{\lie p}$ are   various sets of
semistable points which we study in great detail. In particular,
we have $G$-stable submanifolds $S_\beta$ of $X$ which are called
pre-strata. In case that $\mup$ is proper, the pre-strata form a
decomposition of $X$ and in case that $X$ is compact they are the
strata of a Morse type stratification of $X$. Our results are
generalizations of results of Kirwan obtained in the case that
$G=U\c$ and $X=Z$ is compact.
\end{abstract}

\maketitle

\section{Introduction}

In this paper we continue our study of actions of real reductive
groups on K\"ahler manifolds \cite{HeinznerSchwarz, HSt05}.
Roughly speaking, we extend Kirwan's results on actions of complex
reductive groups (\cite{Kirwan}, see also \cite{Ness} for the
projective case) to the real reductive case. We also obtain new
results in the complex reductive case.

Let $G$ be a closed subgroup of the complex reductive group $H$.
We say that $G$ is real reductive if there is a maximal compact subgroup $U$ of $H$
such that
$K\times \liep\to G$, $(k,\beta)\mapsto k\exp\beta$, is a
diffeomorphism. Here $K:=G\cap U$ is a maximal compact subgroup of $G$ and  $\liep:=\lieg\cap
\im\lieu$ where $\lie u$ denotes the Lie
algebra of   $U$ and $\lieg$ the Lie
algebra of $G$. The Lie algebra $ \lie h$ of  $H$ is the direct sum $\lieu\oplus\im\lieu$. We also say that $G$ is   compatible with
the Cartan decomposition $H=U\c=U\exp \im\lieu$ of $H$. In this paper we fix  $U$ and the real reductive subgroup $G\subset H$.

Assume that $H$ acts holomorphically on a complex K\"ahler manifold
$Z$, that the K\"ahler
form $\omega$ is $U$-invariant and that there is a $U$-equivariant
moment mapping $\mu\colon Z\to\lieu^*$.  For $\xi\in\lieu$ and
$z\in Z$ we set $\mu^\xi(z):=\langle
\mu(z),\xi\rangle:=\mu(z)(\xi)$.  The inclusion
$\im\liep\hookrightarrow \lieu$ induces by restriction a
$K$-equivariant map $\muip\colon Z\to (\im\liep)^*$.  Using a
$U$-invariant inner product on $\lieh$ we can identify
$(\im\liep)^*$ and $\liep$, so we view $\muip$ as a map
$\mu_\liep\colon Z\to\liep$. For $\beta\in\liep$ let
$\mu_\liep^\beta$ denote $\mu^{-\im\beta}$. Then \(
\grad\mup^\beta=\beta_Z \) where $\beta_Z$ is the vector field on
$Z$ corresponding to $\beta$ and $\grad$ is computed with respect
to the Riemannian metric induced by the K\"ahler structure.  We
call $\mup$ the \emph{$G$-gradient map associated with $\mu$}.

For the rest of this paper we fix a $G$-stable locally closed real submanifold $X$ of
$Z$. We may consider $\mu_\liep$ as a mapping $\mup\colon X\to \liep$
such that
\[
\grad\mup^\beta=\beta_X
\]
where the gradient is now computed with respect to the induced
Riemannian metric on $X$. Since $X$ is $G$-stable we have
$\beta_Z(z)=\beta_X(z)$ for $z\in X$. Using the inner product on
$\liep\subset\im\lieu$ we define $\eta_\liep(z):=\frac
12\lVert\mup(z)\rVert^2$, $z\in X$. Let $\cp$ be the set of
critical points of $\eta_\liep$ and $\bp:=\mup(\cp)\subset\liep$.

A strategy for analyzing the $G$-action on $X$ is to view
$\eta_\liep$ as a generalized Morse function in order to obtain a smooth
$G$-stable stratification of $X$ and to study its properties,
as follows.
  Let
$\varphi_t$ denote the flow of the vector field $\grad\eta_\liep$.
For $\beta\in\bp$ consider the set $S_\beta':=\{x\in X\mid
\varphi_t(x)$ has a limit point in
$\cp\cap\mup^{\inv}(K\cdot\beta)$ as $t$ goes to $-\infty \}$. From the Morse theoretical point
of view, the set $S_\beta'$ is a candidate for a Morse-stratum
associated to $\eta_\liep$.  Difficulties arise, since $\cp$ is
almost never smooth, the flow $\varphi_t$ may not exist for
all $t\le 0$ and if it exists a limit point might not be unique.

We get around these difficulties by defining the relevant sets, which we call pre-strata, in terms
of semistability. Set  $\mp:=\mup\inv(0)$ and $\SS_G(\mp):=\{ z\in X: \overline{G\cdot z}\cap
\mup\inv(0)\ne\emptyset\}$, the (open)   set
of semistable points in $X$. If $\beta=0\in\bp$, then we define $S_\beta=S_0=\SS_G(\mp)$. We
define $S_\beta$ for any $\beta\in\bp$ in a similar way, as follows.

For $\beta\in \liep$ let $G^{\beta_+}:=\{g\in G:
\lim_{t\to-\infty}\exp( t\beta)\, g\exp(-t\beta) \text{
exists}\}$. Then $G^{\beta_+}$ is parabolic with Levi component
the centralizer $G^\beta$ of $\beta$. Note that
$G^\beta=K^\beta\exp \liep^\beta$ is compatible with the Cartan
decomposition of $U\c$. Let $X^\beta$ denote the set of zeros of
the vector field $\beta_X$.  This submanifold is stable with
respect to the action of $G^\beta$ and we have a shifted gradient
map $\widehat{\mu_{\liep^\beta}}\colon X^\beta\to \liep^\beta$,
$x\mapsto \mu_{\liep^\beta}(x)-\beta$, whose set of semistable
points is   $\widehat{S_0}:=\{x\in X^\beta:
\overline{G^\beta\cdot
x}\cap\mu_{\liep^\beta}\inv(\beta)\not=\emptyset\}$.  Set
$S^{\beta_+}:=\{z\in X: \lim_{t\to-\infty}\exp(t\beta)\cdot z\in
\widehat{S_0} \}$. Then $S^{\beta_+}$ is $G^{\beta_+}$-stable.  We
call $S_\beta:=G\cdot S^{\beta_+}$ the \emph{pre-stratum associated with
$\beta\in\bp$}.

Let $\twist{G}{G^{\beta_+}}{S^{\beta_+}}$ denote the quotient of
$G\times S^{\beta_+}$ with respect to the $G^{\beta_+}$-action which
is given by $q\cdot (g,s)=(gq\inv,q\cdot s)$ for $g\in G$, $s\in S^{\beta_+}$, $q\in G^{\beta_+}$. In
Section~\ref{section:strata} we obtain our main results concerning
pre-strata.

\begin{introtheorem1} \label{intro:strataslice} The pre-stratum
  $S_\beta$ associated with $\beta$ is a locally closed submanifold of $X$ and the
  natural map
  \[
  \twist G{G^{\beta_+}}S^{\beta_+}\to S_\beta
  \]
  induced by the $G$-action on $X$ is a diffeomorphism.
\end{introtheorem1}

The $G$-action on a pre-stratum has properties similar to those of
$G$-actions on a set of semistable points. The role of the zero fiber
of the gradient map is taken over by the set $\mp(K\cdot
\beta):=\mup\inv(K\cdot \beta)$.  Our second result in
Section~\ref{section:strata} is the following.

\begin{introtheorem2} \label{intro:strataquotient} Let $\beta\in\bp$.
  \begin{enumerate}
  \item The topological Hilbert quotients (see Section \ref{section:slicetheorem})
    $\quot{S^{\beta_+}}{G^\beta}$ and $\quot{S_\beta}{G}$ exist and are
    isomorphic and parameterize the closed orbits.
  \item The inclusions $\m_{\liep^\beta}(\beta)\cap S^{\beta_+}\subset
    S^{\beta_+}$ and $\mp(K\cdot \beta)\cap S_\beta\subset S_\beta$
    induce isomorphisms $(\m_{\liep^\beta}(\beta)\cap
    S^{\beta_+})/K^\beta\cong \quot{S^{\beta_+}}{G^\beta}$ and
    $(\mp(K\cdot \beta)\cap S_\beta)/K\cong \quot{S_\beta}{G}$.
  \item The set of minima of $\eta_\liep\vert S_\beta$ coincides with
    $\mp(K\cdot \beta)\cap S_\beta$.
  \end{enumerate}
\end{introtheorem2}

The content of these results can be summarized in the commutative diagram:

\[
\xymatrix{ & \m_{\liep^\beta}(\beta)\cap S^{\beta_+}
  \ar@{^{(}->}[ld]\ar@{^{(}->}[rr]\ar[dd] \ar@{^{(}->}[rd]& &
  \mp(K\cdot \beta)\cap S_\beta \ar@{^{(}->}[ld]\ar[dd]
  \\
  S^{\beta_+}\ \ar@{^{(}->}[rr] \ar[dd] & & S_\beta \ar[dd] &
  \\
  & (\m_{\liep^\beta}(\beta)\cap S^{\beta_+})/K^\beta \ar@{=}[rr]^\sim
  \ar@{=}[rd]^\sim \ar@{=}[ld]_\sim& & (\mp(K\cdot\beta)\cap S_\beta)/K
  \ar@{=}[ld]_\sim
  \\
  \quot{S^{\beta_+}}{G^\beta}\ar@{=}[rr]^\sim & & \quot{S_\beta}{G} &
}
\]
\\

We also show in Section~\ref{section:strata} that a pre-stratum $S_\beta$
only depends upon $K\cdot \beta$.

\begin{introtheorem3} \label{intro:prestratification} Let $\beta$,
  $\tilde\beta\in\bp$.
  \begin{enumerate}
  \item If $S_{\tilde\beta}\cap \overline{S_\beta}\ne\emptyset$ and
    $S_\beta\not=S_{\tilde\beta}$, then $\|\tilde\beta\|>\|\beta\|$.
  \item The following are equivalent.
    \begin{enumerate}
    \item $K\cdot \beta=K\cdot\tilde\beta$.
    \item $S_{\beta}= S_{\tilde\beta}$.
    \item $S_{\beta}\cap S_{\tilde\beta}\neq\emptyset$.
    \end{enumerate}
  \end{enumerate}
\end{introtheorem3}
If $\mup\colon X\to \liep$ is a proper map, then we show in
Section~\ref{section:proper} that for $\beta\in\bp$, the Morse-theoretical stratum $S_\beta'$ coincides with the pre-stratum $S_\beta$.  Moreover, we
prove the following
\begin{theoremproper}
  If $\mup\colon X\to \liep$ is proper, then $X=\cup_\beta S_\beta$
  where $\beta$ runs through a complete set of representatives of
  $K$-orbits in $\bp$.
  \end{theoremproper}

In Section~\ref{section:compact} we obtain

\begin{theoremcompact} If $X$ is compact, then there are only finitely many pre-strata.
\end{theoremcompact}

If $X=Z$ is compact and $G=U\c$, then most of the above is proved in
\cite{Kirwan}. This holomorphic and compact case has been studied
quite intensively. For example, the strata $S_\beta$ are smooth analytic
subsets of $Z$ which are open in their closure with respect to the analytic Zariski topology,  there is only
one open pre-stratum $S_{\beta_0}$ (if $X$ is connected) and its complement is a closed analytic
subset of $X$.  Also, the open pre-stratum $S_{\beta_0}$ is connected and is
minimal in the sense that $\beta_0$
is a closest point
 to zero
in $\mup(X)\subset\liep$.
Many of these results about open pre-strata hold more generally
for meromorphic actions of complex reductive groups on compact K\"ahler
manifolds (see \cite{Sommese, Fujiki}).

An important application noted in \cite{Kirwan} are the Morse equalities in
equivariant cohomology in the holomorphic compact case.
Morse inequalities also hold in our setup for general real reductive group actions
on compact manifolds. We illustrate in Section~\ref{section:morse} by simple
examples that equalities are almost never valid for actions of real reductive
groups.

The simplest interesting situation where our results apply and is not co\-vered by
previously known results is the real projective case. Here
$X\subset\PP(V_\RR)$ where $V_\RR$ is a representation space of an algebraic
semisimple real group $G$ and $X$ is $G$-stable and closed.  In this
setup we may choose
$Z=\PP(V)$ where $V=V_\RR\otimes \CC$  and $\PP(V)$ is given the Fubini-Study
metric.

The special case $X=\PP(V_\RR)$ and $G$ real semisimple is treated
in \cite{Marian}. However, there is an error in the proof of
Theorem 1: the assumption on the bottom of page 783 that one may
assume that $k_1=k_2=e$ is incorrect.

Still in the projective setup we have the following result of Sch\"utzdeller
\cite{Schuetzdeller} which is new even in the case that $G=U\c$ and which
generalizes results of Kostant (\cite{Kos73}), Atiyah (\cite{Ati82}),
Guillemin and Sternberg (\cite{GS82}), Mumford (Appendix of \cite{Ness}),
Kirwan (\cite{Kir84b}) and O'Shea and Sjamaar (\cite{OS00}).
\\

\textit{Let $\liea_+$ be a positive Weyl chamber for the $K$-action on
  $\liep$. Let $Y$ be a semi-algebraic closed $G$-stable subset of
  $\PP(V)$ which is irreducible (in a suitable sense). Then
  $\mup(Y)\cap\liea_+$ is a convex polytope.}
\\

The proof of this convexity result makes essential use of the
stratification obtained here.

Another related interesting simple case is where $X=Z$ is
$U$-homogeneous and $G$ is a real form of $U\c$.  Our results applied
to $G$ and also to the complex reductive group $K\c$ imply
Matsuki duality (\cite{Ma,Vilonen,Bremigan}, see Corollary \ref{corollary:matsuki}) and may be viewed as a generalization thereof.

Most of our results in the holomorphic non-compact case are new.  In this
setup $G$ is a complex reductive group and $X$ is a K\"ahler manifold.
The re\-le\-vant maps are holomorphic
and the quotients discussed above are K\"ahler spaces. In particular,
a pre-stratum $S_\beta$ is
$G$-equivariantly and biholomorphically identified with a holomorphic
$G$-bundle over the compact complex manifold $G/G^{\beta_+}$.

Also, in the general real case, a pre-stratum $S_\beta$ fibers
over the $G$-ho\-mo\-ge\-ne\-ous compact manifold $G/G^{\beta_+}$
which now is a $K$-homogeneous Riemannian manifold. The open pre-strata are of
special interest. In many
situations, including the case where $\mup$ has closed image,
there is a minimal pre-stratum $S_{\beta_0}$
in the sense that $\beta_0$
is a closest point to zero in
$\mup(X)$.
This pre-stratum is always open and usually serves as a good substitute for the set
of semistable points which quite often is empty.

In principle there could be many different minimal pre-strata and also
other open pre-strata which are not minimal in the above sense. We do not
know any example of a connected manifold $X$ where this is indeed the
case. Any such example would imply that $\mup(X)\cap \liea_+$ is not
a convex subset of $\liea_+$. In situations where the convexity result
of Sch\"utzdeller applies, the
minimal pre-stratum is unique. The uniqueness of the minimal pre-stratum holds
in complete generality if $X=Z$ and $G=U^\ce$. In this case very
general convexity theorems are available
(\cite{HHuck}).

Another interesting extreme situation occurs when there is a
maximal pre-stratum $S_\beta$. By this we mean that
$\beta\in\mu_\liep(X)$ and $\norm\beta\ge\norm{\mup(z)}$ for all
$z\in X$. In this case the function $\eta_\liep$ is constant on
$S_\beta$ and we explain in Section~\ref{section:applications}
that the general theory implies that $S_\beta=\mp(K\cdot
\beta)=\mup\inv(K\cdot \beta)$ and that
$S^{\beta_+}=\mp(\beta)=\mup\inv(\beta)$. Every $G$-orbit in
$S_\beta$ is a $K$-orbit and every $G^\beta$-orbit in
$S^{\beta_+}$ is a $K^\beta$-orbit. Moreover,
$S^{\beta_+}=\mp(\beta)$ is a union of connected components of
$X^\beta$ and the Slice Theorem gives a $K$-equivariant
diffeomorphism $\twist{K}{K^\beta}{\mp(\beta)}\cong S_\beta$.
These considerations generalize an old theorem of Wolf
(\cite{Wolf}, see Corollary \ref{corollary:old wolf}) which says
that a real form $G$ of a complex semisimple group acting on a
generalized flag manifold has a unique closed $G$-obit.

\medskip
\noindent {\sc Acknowledgements}
\medskip

\noindent The authors thank the referee for helpful remarks.

\section{Compatible subgroups and gradient maps}

Let $U$ be a compact Lie group. Then $U$ has a natural real linear
algebraic group structure, and we denote by $U\c$ the
corresponding complex linear algebraic group (\cite{Chevalley}).
The group $U\c$ is reductive and is the universal complexification
of $U$ in the sense of \cite{Hochschild}.  On the Lie algebra
level we have the Cartan decomposition $\lieu\c=\lieu +\im\lieu$
with a corresponding Cartan involution $\theta:\lieu\c\to\lieu\c$,
$\xi+\im\eta\mapsto \xi - \im\eta$, $\xi$, $\eta\in\lieu$. The
real analytic map $U\times \im\lieu\to U\c, \quad (u,\xi)\mapsto
u\exp \xi$, is a diffeomorphism.  We refer to the decomposition
$U\c=U\exp(\im\lieu)\cong U\times \im\lieu$ as the \emph{Cartan
decomposition of $U\c$}, and we fix it for the remainder of this
paper.

Let $G$ be a real Lie subgroup of $U\c$. We say that $G$
\emph{is a compatible subgroup of} $U^\ce$ if it is compatible
with the Cartan decomposition of $U\c$, i.e., if
\[
K\times \liep\to G,\qquad (k,\beta)\mapsto k\exp\beta,
\]
 is a diffeomorphism where $K=G\cap U$ and
 $\liep=\lieg\cap\im\lieu$.
In particular, $G$ is a closed subgroup of $U\c$ if and only if $K$ is
compact.

Let $Z$ be a smooth holomorphic $U\c$-space, i.e., $Z$ is a
complex manifold with a holomorphic action $U\c\times Z\to Z$. We
assume that $Z$ has a K\"ahler form $\omega$ which is
$U$-invariant. We also assume that we are given a $U$-equivariant
moment mapping $\mu\colon Z\to\lieu^*$.  By definition we then
have
\[
d\mu^\xi=\imath_{\xi_Z}\omega
\]
where $\mu^\xi$ is $\mu$ followed by evaluation at $\xi\in\lieu$ and
$\xi_Z$ is the vector field on $Z$ with one parameter group
$(t,z)\mapsto \exp t\xi\cdot z$.

Now let $G$ be a compatible closed subgroup of $U\c$. We have the
subspace $\im\liep\subset\lieu$ and a corresponding mapping
$\muip\colon Z\to(\im\liep)^*$ where $\muip$ is just $\mu$
followed by restriction to $\im\liep$.  The map $\muip$ is the
correct analogue of $\mu$ when one is considering the action of
$G$ rather than the action of $U\c$. In order to simplify notation
we replace consideration of $\muip$ by that of $\mup\colon Z\to
\liep$ where
\[
\mup^\beta(z):=\langle \mup(z),\beta \rangle:=\muip^{-\im\beta}(z).
\]
Here $\langle\cdot,\cdot\rangle$ denotes a $K$-invariant inner product on
$\liep\subset\im\lieu$. The map $\mup$ is called the \emph{ $G$-gradient map associated with
  $\mu$}  since the
   equation $d\mu^\xi=\imath_{\xi_Z}\omega$ is equivalent to \( \grad \mup^\beta=\beta_Z
\) for all $\beta\in \liep$. Here the gradient is computed with
respect to the Riemannian metric $\left(.,.\right)$ on $Z$ given by
$\left(v,w\right)=\omega(v,Jw)$ for all $z\in Z$ and $v,w\in T_z(Z)$
where $J$ denotes the complex structure on $T(Z)$.

For the rest of this paper we fix a $G$-stable locally closed
submanifold $X$ of $Z$. Of course, a special important case is
where $X=Z$. From now on we also denote  the restriction of $\mup$
to $X$ by $\mup$.   We have
\[\grad \mup^\beta=\beta_X\]
where $\grad$ is computed with respect to the induced Riemannian
structure.  Similarly, $\perp$ denotes perpendicularity relative to
the induced Riemannian metric on $X$. We have the following two
elementary results (see \cite[Lemmas 5.1, 5.4] {HeinznerSchwarz}).

\begin{lemma} \label{lemma:monotonic} Let $x\in X$ and $\beta\in
  \liep$.  Then either $\beta_X(x)=0$ or the function $t\to
  \mup^\beta(\exp t\beta\cdot x)$ is strictly increasing.
\end{lemma}

For a subspace $\liem$ of $\lieu^\ce$ and $x\in X$ let $\liem\cdot x$  denote the
subspace $\{\beta_X(x): \beta\in \liem\}$ of the tangent space
$T_x(X)$.

\begin{lemma} \label{lemma:momentumkern} For all $x\in X$ we have
  $\ker d\mup(x)=(\liep\cdot x)^\perp$.
\end{lemma}

We will use $\lVert\cdot \rVert$ to denote the norm functions associated to
$\langle\cdot,\cdot\rangle$ and $(\cdot,\cdot)$. The critical points
of the norm square function $\eta_{\liep}(x):=\frac
12\lVert\mup(x)\rVert^2$ will be of central importance in the rest of
this paper.

\begin{lemma}\label{lemma:gradients}
  Let $x\in X$ and $\beta=\mup(x)$. Then
  $\grad\eta_{\liep}(x)=\beta_X(x)$.
\end{lemma}
\begin{proof} This follows from $
  d\eta_{\liep}(x)(v)=\langle\beta,d\mup(x)v\rangle=d\mup^\beta(x)(v)=(v,\beta_X(x))
  $.
\end{proof}

\begin{corollary} \label{corollary:critical} Let $x\in X$ and set
  $\beta:=\mup(x)$.  The following are equivalent.
  \begin{enumerate}
  \item $\beta_X(x)=0$,
  \item $d\mu_\liep^\beta(x)=0$ and
  \item $d\eta_{\liep}(x)=0$.
  \end{enumerate}
\end{corollary}

Let $x\in X$ and let $\beta\in\liep_x:=\{\beta\in\liep:
\beta_X(x)=0\}$. Differentiating the action of $G_x$ on $T_x(X)$
gives rise to a linear action of $\lieg_x$ on $T_x(X)$. Since
elements of $i\lieu_x$ act as selfadjoint operators on $T_x(Z)$
relative to the induced Riemannian structure, our element $\beta$
acts on $T_x(X)$ as a selfadjoint operator $d\beta_X(x)$ and has
real eigenvalues.
\begin{proposition} \label{proposition:hessian} Let $v\in T_x(X)$ be
  an eigenvector of $\beta\in \liep_x$ with eigenvalue
  $\lambda(\beta)$.  Let $\gamma(t)$ be a smooth curve in $X$ with
  $\gamma(0)=x$ and $\frac d{dt}\gamma(0)=v$. Then
  \begin{enumerate}
  \item
    $\frac{d^2}{dt^2}\,(\mup^{\beta}\circ\gamma)(0)=\lambda(\beta)\lVert
    v\rVert^2$.\label{item:propostion:hessian1}
  \item If $x$ is a critical point of $\eta_{\liep}$ and
    $\beta:=\mup(x)$, then
    \[\frac{d^2}{dt^2}\,(\eta_\liep\circ\gamma )(0)=
    \lambda(\beta)\lVert v\rVert^2+\lVert d\mu_{\liep}(x)(v)\rVert^2.
    \]\label{item:propostion:hessian2}
  \end{enumerate}
\end{proposition}

\begin{proof} We have
  \[
  \frac d{dt}\mup^\beta(\gamma(t))=\left(\frac
  d{dt}\gamma(t),\beta_X(\gamma(t))\right)_{\gamma(t)}
  \]
  where $(\cdot,\cdot)_z$ denotes the inner product on $T_z(X)$ at
  $z\in X$. There is a  neighborhood  $U$ of $x\in X$ and a local diffeomorphism $\Psi\colon U\to T_x(X)$ where $\Psi(x)=0$ and
  $d\Psi(x)=\Id$.  Using the local coordinates given by $\Psi$
  we have $\left(\ ,\ \right)_{\gamma(t)}=\left(\ ,\ \right)_x+t\left(\ ,\ \right)_t$
  where $\left(\ ,\ \right)_t$ is a bilinear form on $T_x(X)$ depending smoothly on $t$.
  From $\beta_X(\gamma(t))=t\cdot d\beta_X(x)\cdot (\frac{d}{dt}\gamma)(0)+t^2R_0(t)
  =t\lambda(\beta)v+ t^2R_0(t)$, where $R_0(t)$ is smooth, we get
\[
\begin{split}
\frac d{dt}(\mup^\beta\circ\gamma)(t)
& =\left(\frac d{dt}\gamma(t),\beta_X(\gamma(t))\right)_x+
t\left(\frac d{dt}\gamma(t),\beta_X(\gamma(t))\right)_t \\
& =\left(\frac d{dt}\gamma(t), t\lambda(\beta)v +t^2R_0(t)\right)_x
+t\left(\frac d{dt}\gamma(t), t\lambda(\beta)v+t^2R_0(t)\right)_t  \\
& =\left(\frac d{dt}\gamma(t),t\lambda(\beta)v\right)_x+t^2R_1(t).
\end{split}
\]
This implies (\ref{item:propostion:hessian1}).

If $\beta$ is as in (\ref{item:propostion:hessian2}), then $\mup(\gamma(t))=\beta+t\cdot
d\mup(x)\cdot v+t^2R_2(t)$. Hence
\[
\begin{split}
\frac d{dt}(\eta_{\liep}\circ\gamma)(t) & =
\left<d\mup(\gamma(t))\frac d{dt}\gamma(t), \mup(\gamma(t))\right>\\
& =
\frac d{dt}(\mup^\beta\circ\gamma)(t) +t\left< d\mup(\gamma(t)) \frac d{dt}\gamma(t), d \mup(x)\cdot v\right>+t^2R_3(t)
\end{split}
\]
and   (\ref{item:propostion:hessian2}) follows.
\end{proof}

If $\liem$ is a subspace of a Lie algebra $\lie l$ and
$\beta\in\lie l$, set $\liem^\beta:=\{\xi\in \liem:
[\xi,\beta]=0\}$.

\begin{lemma}\label{lemma:hessian2}
  Let $\beta:=\mu_\liep(x)$ and assume that $\beta\in\liep_x$. Let $\zeta\in\liek$. If  $\zeta_X(x)$
  is the sum of positive eigenvectors of $\beta$, then
  $\zeta\in\liek^\beta$.
\end{lemma}

\begin{proof}
  Let $\gamma(t):=\exp t\zeta\cdot x$.
  Proposition~\ref{proposition:hessian} implies  that \( 0\le
  \tfrac{d^2}{dt^2}(\mu_\liep^\beta\circ \gamma)(0) \).  Since $\mup$
  is $K$-equivariant and the $K$-action on $\liep$ is by linear
  isometries we have
  \[
  \tfrac{d^2}{dt^2}(\mu_\liep^\beta\circ\gamma)(0) = \tfrac{d^2}{dt^2}\vert_0 \langle \exp(t\zeta)\cdot \mu_{\liep}(x),\beta\rangle =\langle
  [\zeta,[\zeta,\beta]],\beta\rangle=-\lVert [\zeta,\beta]\rVert^2\le
  0.
  \]
  This shows that $[\zeta,\beta]=0$, i.e., $\zeta\in\liek^\beta$.
\end{proof}

\section{The Slice and Quotient Theorems} \label{section:slicetheorem}

In this section we recall for the convenience of the reader some
results from \cite{HeinznerSchwarz,HSt05}. For any Lie group $G$,
a closed subgroup $H$ and any set $S$ with an $H$-action we denote
by $\twist GHS$ the $G$-bundle over $G/H$ associated with the
$H$-principal bundle $G\to G/H$. This is the orbit space of the
$H$-action on $G\times S$ given by $h\cdot (g,s)=(gh\inv,h\cdot
s)$ where $g\in G$, $s\in S$ and $h\in H$. The $H$-orbit of $(g,s)$, considered as a point in $\twist
GHS$, will be denoted by $[g,s]$.

Let $G=K\exp\liep$ be a compatible closed subgroup of $U\c$.  For
$\beta\in\liep$ set $\mp(\beta):=\mup\inv(\beta)\subset X$ and set
$\mp:=\mp(0)$.  Let $x\in \mp$.  Then $G_x=K_x\exp\liep_x$
(\cite[5.5]{HeinznerSchwarz}). Since the $G_x$-representation on
$T_x(X)$ is completely reducible (\cite[14.9]{HeinznerSchwarz}),
there is a $G_x$-stable decomposition $T_x(X)=\lieg\cdot x\oplus
W$. Now the Slice Theorem for $Z$ (\cite[14.10,
14.21]{HeinznerSchwarz}) pulls back to the following Slice Theorem
for $X$:

\begin{theorem}[(Slice Theorem)] \label{slicetheorem} Let $x\in \mp$.
  Then there exists a $G_x$-stable open neighborhood $S$ of $0\in W$, a
  $G$-stable open neighborhood $\Omega$ of $x\in X$ and a
  $G$-equivariant diffeomorphism $\Psi\colon \twist G{G_x}S\to \Omega$
  where $\Psi([e,0])=x$.
\end{theorem}

Actually, we have a Slice Theorem at every $x\in X$. Set
$\beta:=\mup(x)$ and let $G^\beta=\{g\in G: \Ad g\cdot
\beta=\beta\}$ denote the centralizer of $\beta$. Then we have a slice
for the action of $G^\beta$, as follows.

The centralizer $G^\beta$ is a compatible subgroup of $U\c$ with
Cartan decomposition $G^\beta=K^\beta\exp(\liep^\beta)$ where
$K^\beta=K\cap G^\beta$ and $\liep^\beta=\{\xi\in\liep:
\ad(\xi)\beta=0\}$. The group $G^\beta$ is also compatible with
the Cartan decomposition of $(U\c)^\beta=(U^\beta)\c$ and $\beta$
is fixed by the action of $U^\beta$ on $\lieu^\beta$. This implies
that the $\lieu^\beta$-component of $\mu$ defines a
$U^\beta$-equivariant shifted gradient map
$\widehat{\mu_{\lieu^\beta}}\colon Z\to \lieu^\beta$,
$\widehat{\mu_{\lieu^\beta}}(z)=\mu_{\lieu^\beta}(z)-\beta$. The
associated $G^\beta$-gradient map is given by
$\widehat{\mu_{\liep^\beta}}\colon X\to \liep^\beta$,
$\widehat{\mu_{\liep^\beta}}(z)=\mu_{\liep^\beta}(z)-\beta$. This
shows that the Slice Theorem applies to the action of $G^\beta$ at
every point $x\in
\widehat{(\mu_{\liep^\beta}})\inv(0)=\mathcal{M}_{\liep^\beta}(\beta)$.
In particular, if $G$ is commutative, then we have a Slice Theorem
for $G$ at every point of $X$.

For $H$ a subgroup of $G$, $M$ a subset of $Z$ and $Y$ an
$H$-stable subset of $Z$, we define the saturation $\SS_H(M)(Y)$
to be $\{z\in Y: \overline{H\cdot z}\cap M\neq\emptyset\}$.  Here
$\overline{H\cdot z}$ denotes the closure of $H\cdot z$ in $Y$. In
general, $\SS_H(M)(Y)$ is a proper subset of $\SS_H(M)(Z)\cap Y$.
But in the case where $Y$ is closed in $Z$ these sets agree. The
set $\SS_G(\mp)(Y)$ is called the \emph{set of semistable points}
of $Y$ with respect to $\mup$. In (\cite{HSt05}) it has been shown
that $\SS_G(\mp(\beta))(Z)$ is open in $Z$ for every
$\beta\in\liep$. For a closed $G$-stable subset $Y$ of $Z$ this
implies that $\SS_G(\mp(\beta))(Y)$ is open in $Y$. Inspecting the
proof in \cite{HSt05} one can more generally show that
$\SS_G(\mp(\beta))(X)$ is open in $X$ for any locally closed
$G$-stable submanifold $X$ of $Z$. Since we have fixed $X$, we
usually write $\SS_G(\mp(\beta))$ for $\SS_G(\mp(\beta))(X)$. The
set of semistable points plays a major role in
\cite{HeinznerSchwarz}. One reason for this is the quotient
theorem which we now formulate.

Let $Y\subset Z$ be $G$-stable and let $x$, $y\in Y$. We define a
relation $\sim$ on $Y$ where $x\sim y$ if and only if $Y\cap
\overline{G\cdot x}\cap\overline{G\cdot
  y}\neq\emptyset$.
  If this relation
is in fact an equivalence relation we denote the corresponding
quotient by $\quot YG$ and call it the \emph{topological Hilbert quotient of $Y$ by the action of $G$}.

\begin{theorem}[(Quotient Theorem \cite{HeinznerSchwarz})]
  \label{theorem:quotient} Assume that $X=\SS_G(\mp)$. Then the
  topological Hilbert quotient $\quot XG$ exists and has the following
  properties.
  \begin{enumerate}
  \item Every fiber of $\pi$ contains a unique closed $G$-orbit.  Any other orbit in the fiber has
    strictly larger dimension.
  \item The closure of every $G$-orbit in a fiber of $\pi$ contains the
    closed $G$-orbit.
  \item Every fiber of $\pi$ intersects $\mp$ in a unique $K$-orbit
    which lies in the unique closed $G$-orbit.
  \item The inclusion $\mp\hookrightarrow X$ induces a homeomorphism
    $\mp/K\cong \quot XG$.
  \end{enumerate}
\end{theorem}

As in the case of the Slice Theorem, we have local versions of the
Quotient Theorem. For any $\beta\in \liep$ we have the open subset
of semistable points $\SS_{G^\beta}(\m_{\liep^\beta}(\beta))$ in
$X$, and we can apply the Quotient Theorem for the action of
$G^\beta$.

\section{Fixed points and parabolic
  subgroups}\label{section:parabolic}

Let $\beta\in\liep$. We have a vector field $\beta_G$ whose
one-parameter subgroup is given by\\
$(t,y)\mapsto
\exp(t\beta)y\exp(-t\beta)$. Then
$$G^\beta=\{y\in G: \beta_G(y)=0\}=\{y\in G: \exp(t\beta)y\exp(-t\beta)=y \text{ for
  all $t\in \RR$}\}.
  $$
We have the  parabolic subgroup
\[
G^{\beta_+}:=\{y\in G:
\lim_{t\to-\infty}\exp(t\beta)y\exp(-t\beta)\text{ exists}\}
\]
with unipotent radical
\[
R^{\beta_+}:=\{y\in G:
\lim_{t\to-\infty}\exp(t\beta)y\exp(-t\beta)=e\}.
\]
Then $G^{\beta_+}$ is the semi-direct product of $G^\beta$ with
$R^{\beta_+}$ and we have the projection
  $\pi^{\beta_+}\colon G^{\beta_+}\to G^\beta$,
$\pi^{\beta_+}(y):=\lim_{t\to-\infty}\exp(t\beta)y\exp(-t\beta)$.
\begin{lemma}\label{lemma:GgleichKGbeta}
  For every $\beta\in\liep$ we have $G=KG^{\beta_+}$.
\end{lemma}

\begin{proof}
  The adjoint orbit $O$ of $U$ through $\beta\in\im\lieu\cong \lieu$
  can be considered as a K\"ahler
  manifold endowed with a
  holomorphic $U\c$-action. In particular, $G$-acts on $O$.  In
  \cite{HSt05} it is shown that $G\cdot\beta=K\cdot\beta$.
 Since $G\cdot\beta\cong G/G^{\beta_+}$ this proves
  the claim.
\end{proof}

We introduce submanifolds of $X$ analogous to $G^\beta$ and
$G^{\beta_+}$. For $\beta\in\liep$ we have the corresponding
vector field $\beta_X$ on $X$ whose one-parameter subgroup is
given by $(t,y)\mapsto \exp(t\beta)\cdot y$.  We have the set
\[
X^\beta:=\{y\in X: \beta_X(y)=0\}
\]
which is the set of fixed points $\{y\in X: \exp t\beta\cdot
y=y\text{ for all $t\in \RR$}\}$. The set $X^\beta$ is
$G^\beta$-stable and is a subset of the $G^{\beta_+}$-stable set
\[
X^{\beta_+}:=\{y\in X: \lim_{t\to - \infty}\exp t\beta\cdot y\text{
  exists}\}.
\]
The map $p^{\beta_+}\colon X^{\beta_+}\to X^\beta$,
$p^{\beta_+}(y)=\lim_{t\to -\infty}\exp t\beta\cdot y$ is well
defined, $G^\beta$-equivariant, surjective and its fibers are
$R^{\beta_+}$-stable.
\\

In the generality above, the map $p^{\beta_+}$ may not be continuous.
\begin{example}
Let $X=\PP^1(\CC)$ with the $\RR$-action given by  $t\cdot [z,w]=[tz,t\inv w]$, $t\in\RR$, $[z,w]\in \PP^1(\CC)$.
The infinitesimal generator of this action is $\beta=\left(\begin{smallmatrix} 1 & 0 \\ 0 & -1\end{smallmatrix}\right)$. Then $X^{\beta_+}=X$, $X^{\beta}=\{[1,0], [0,1]\}$, $p^{\beta_+}([z,w])=[1,0]$ when $z\neq 0$ and $p^{\beta_+}([0,1])=[0,1]$. So $p^{\beta_+}$ is not continuous. Moreover, the Hilbert quotient $\quot X\RR$ does not exist. Everything is fine, however,  if we restrict to a set of semistable points as we do below.
\end{example} 

In the following we fix   $\beta\in\liep$, $\|\beta\|=1$ and
discuss, in the spirit of this paper, the relevant properties of
the action of the group $\Gamma:=\exp\gamma$ on $X$ where
$\gamma:=\RR\beta$. Note that $\Gamma$ is a closed compatible
subgroup of $U^\CC$. As an abstract Lie group $\Gamma$ is just the
additive group $\RR$. The $\Gamma$-gradient map on $X$ is given by
$\mu_\gamma(y)=\mu_\gamma^\beta(y)\beta$ and will in the following
be identified with $\mu_\gamma^\beta\colon X\to \RR$.   The
isotropy group $\Gamma_x$ of $\Gamma$ at every point $x\in X$ is
compatible with the Cartan decomposition of $U\c$. Since
$\Gamma\cong \RR$,
we have that either $\Gamma_x=\{0\}$ or
$\Gamma_x=\Gamma$. 
Since $\Gamma$ is commutative, there is a slice for the $\Gamma$-action at every point of $Z$.
The Slice Theorem applied at $x\in X\setminus
X^\beta$ gives the existence of an open $\Gamma$-stable
neighborhood $\Omega$ of $x$, a closed submanifold $S$ of $\Omega$
with $x\in S$ such that the map $\Psi\colon \Gamma\times S\to
\Omega$, $(g,s)\mapsto g\cdot s$ is a diffeomorphism.

In the case that $x\in X^\beta$, the Slice Theorem gives a
linearization of the $\Gamma$-action near $x$, as follows.  The
linearized vector field $d\beta_X(x)$ acts on $T_x(X)$ as a
self-adjoint operator, also denoted by $\beta$, and $W:=T_x(X)$ is
a direct sum of one dimensional eigenspaces. We define
$W^{\beta_+}:=W^\beta\oplus W^+$ where $W^+$ is the sum of the
eigenspaces of $\beta$ with positive eigenvalues and $W^\beta$ is
the zero eigenspace of $\beta$. By the Slice Theorem there is an
open $\Gamma$-stable neighborhood $S$ of zero in $T_x(X)$ and a
$\Gamma$-equivariant diffeomorphism $\Psi$ which maps $S$ onto an
open  neighborhood $\Omega$ of $x$ in $X$. Since $S$ is
$\Gamma$-stable, it  contains $(S\cap W^\beta)\oplus W^+$. 

Suppose that $X$ coincides with a set of semistable points, e.g., suppose that for some $r\in\RR$
$$
X=\{z\in X\colon \overline{\Gamma\cdot z}\cap(\mu^\beta_\gamma)\inv(r)\}\neq\emptyset\}.
$$
Then the topological Hilbert quotient $\pi\colon X\to\quot X\Gamma$ exists and we may choose $S$ such that $\Psi(S)=
\Omega\subset X$ is saturated with respect to $\pi$. Then $\Psi$ identifies $S\cap W^{\beta_+}=(S\cap W^\beta)\times W^+$ with $\Omega\cap X^{\beta_+}$ and 
 $S\cap W^\beta$   with $\Omega\cap X^\beta$.
Moreover, $(p^{\beta_+})\inv\Psi(S\cap W^\beta) = \Psi((S\cap
W^\beta) \oplus W^+)$ is   closed in  $\Omega$.

\begin{remark}
  The Slice Theorem applied to $\Gamma$ at $x\in X^\beta$ shows that, for $y$ near $x$,
  the limit
  $\exp t_n\beta\cdot y$ exists for some
  sequence $t_n\to -\infty$ if and only if $\lim_{t\to-\infty} \exp
  t\beta\cdot y=x$.
\end{remark}

We summarize our discussion as follows.

\begin{proposition} \label{proposition:strata} 
Let $\beta\in\liep$ and 
  $\Gamma=\exp\RR\beta$.
Let $r\in \RR$ and set $X_r:=\{z\in X\colon \overline{\Gamma \cdot z}\cap(\mu^\beta_\gamma)\inv(r)\neq\emptyset\}$ and set $X^\beta_r=X^\beta\cap X_r$. Then $X_r$ is open in $X$ and we assume that it is nonempty. Set
$$
X^{\beta_+}_r:=\{z\in X_r\colon\lim_{t\to-\infty}\exp(t\beta)\cdot z\in  X_r^\beta\}.
$$
  \begin{enumerate}
  \item The set $X^\beta$ is a closed submanifold of $X$.
  \item The set $X_r^{\beta_+}$ is a locally closed $G^{\beta_+}$-stable submanifold of
    $X$
  \item The tangent space of $X_r^{\beta_+}$ at $x\in X_r^\beta$ is
   $W^{\beta_+}$.\label{item:proposition:strata3}
  \item The map $p^{\beta_+}$ is a $G^\beta$-equivariant strong
    deformation retraction of $X^{\beta_+}_r$ onto $X_r^\beta$.
  \item $p^{\beta_+}$ is $G^{\beta_+}$-equivariant where the action of
    $G^{\beta_+}$ on $ X_r^\beta$ is via the quotient morphism to $G^\beta$.
   \end{enumerate}
\end{proposition}
\begin{remark}
  In addition to the properties in
  Proposition~\ref{proposition:strata}, we have that the map
  $p^{\beta_+}$ is a smooth locally trivial fibration which realizes
  $X_r^{\beta}$ as the topological Hilbert quotient of $X_r^{\beta_+}$
  with respect to the action of $\Gamma$.  Each fiber
  $F^{\beta_+}_x:=(p^{\beta_+})\inv(x)$ of $p^{\beta_+}$ is
  $(G^\beta)_x$-equivariantly diffeomorphic to a
  $(G^\beta)_x$-representation $W^+(x)$ which, up to isomorphism,
  depends only on the connected component of $x$ in $X^\beta\cap X_r$.
 Furthermore, the map $\twist
  {G^\beta}{(G^\beta)_x}F_x^{\beta_+}\to
  (p^{\beta_+})\inv(G^\beta\cdot x)$, $[a,z]\mapsto a\cdot z$ is a
  $G^\beta$-equivariant diffeomorphism.
\end{remark}

\begin{remark}
  If $G=K\c$ for some compact subgroup $K$ of $U$ and if $X$ is a complex
  submanifold of $Z$, then the manifolds discussed in
  Proposition~\ref{proposition:strata} are complex analytic and the
  maps are holomorphic.
\end{remark}


The results above show that $X^\beta_r$ is open in $X^\beta$. It is also closed.

\begin{proposition} \label{proposition:strataminimum} Let
  $\beta\in\liep$.
  \begin{enumerate}
  \item The function 
  $\mu^\beta_\gamma$ is locally
    constant on $X^\beta$. In particular, $X^\beta_r$ is open and closed in
    $X^\beta$. \label{item:proposition:strataminimum1}
  \item We have 
  $\mu^\beta_\gamma(z)\ge r$ for $z\in X_r^{\beta_+}$ and
    equality holds if and only if $z\in X^\beta_r$.\label{item:proposition:strataminimum2}
  \end{enumerate}
\end{proposition}

\proof For $v\in T_y(X)$ and $z\in X^\beta$,
$d\mu^\beta_\gamma(z)=\left(\beta_X(z), v\right)=0$. 
Hence 
$\mu^\beta_\gamma$ is locally constant on $X^\beta$ and we have (\ref{item:proposition:strataminimum1}).
Let $z\in X_r^{\beta_+}$. Then 
$\mu^\beta_\gamma(\exp(t\beta)\cdot z\ge\lim_{t\to -\infty}\mu^\beta_\gamma(\exp t\beta\cdot y)=r$. 
Note that
$\mu^\beta_\gamma(z)=r$ if and only if $\exp t\beta\cdot z=z$ for all $t\in
\RR$ (Lemma~\ref{lemma:monotonic}), i.e., $z\in X^\beta$. This
establishes (\ref{item:proposition:strataminimum2}). \qed

\begin{remark} By replacing $\beta$ with $-\beta$ one obtains
  analogous results for\\ $X_r^{\beta_-}:=\{z\in X_r: \lim_{t\to+\infty}\exp
  t\beta\cdot z\text{ exists and lies in }X^\beta_r\}$.
\end{remark}

Let $\beta\in\liep$ and let $\lieg^\beta=\liek^\beta+\liep^\beta$
denote the Cartan decomposition of the Lie algebra $\lieg^\beta$ of
the centralizer $G^\beta$. For $r\in \RR$ let
$H_r(\beta):=\{\zeta\in\liep^\beta: \langle\beta,\zeta\rangle=r\}$. In
the following it is important to note that $H_r(\beta)$ is a
hyperplane in $\liep^\beta$ if and only if $\beta\not=0$. If
$\beta\not=0$ we have the corresponding half-spaces
$H^+_r(\beta):=\{\zeta\in\liep^\beta: \langle\beta,\zeta\rangle\ge
r\}$ and $H^-_r(\beta):=\{\zeta\in\liep^\beta:
\langle\beta,\zeta\rangle\le r\}$. These sets are also defined in the
case where $\beta=0$ but then they are either empty ($r\neq 0$) or all
of them coincide with $\liep^\beta$ ($r=0$).

\begin{corollary} \label{corollary:hyperplane} Let $\beta\in\liep$.
  \begin{enumerate}
  \item $\mu_{\liep^\beta}(X_r^{\beta_+})\subset H^+_r(\beta)$,
  \item $\mu_{\liep^\beta}(X_r^{\beta_-})\subset H^-_r(\beta)$ and
  \item $X_r^\beta=X_r^{\beta_-}\cap X_r^{\beta_+}=
    \mu_{\liep^\beta}\inv(H_r(\beta))\cap X_r^{\beta_+}$.
  \end{enumerate}
\end{corollary}

\begin{remark} \label{remark:independent} For fixed $r$ the manifold
  $X_r^{\beta_+}$ is defined only in terms of the group
  $\Gamma:=\exp\RR\beta$.  This means that $X_r^{\beta_+}$ remains
  unchanged for any compatible subgroup $G$ of $U^\CC$ which contains
  $\Gamma$ and stabilizes $X$. 
\end{remark}

\begin{remark}
The set $X_r^\beta$ is $G^\beta$-stable, which one sees as follows. The map $\mu_{\liep^\beta}\colon X^\beta\to\liep^\beta$ is $K^\beta$-equivariant where $K^\beta$ is the centralizer of $\beta$ in $K$. Thus $\mu^\beta\colon X^\beta\to \RR\cdot\beta$ is also $K^\beta$-equivariant. Since $G^\beta=K^\beta\exp(\liep^\beta)$, $K^\beta$ maps onto the component group of $G^\beta$, hence $X^\beta_r$ is $G^\beta$-invariant for any $r\in \RR$.
\end{remark}

\section{Critical points, Slices and Quotients}
\label{section:strata}

Let $\mu_{\liep^\beta}\colon X\to\liep^\beta$ denote the gradient
map associated with the moment map $\mu_{\lieu^\beta}\colon Z\to
\lieu^\beta$ (see Section~\ref{section:slicetheorem}). Recall that
$\m_{\liep^\beta}(\beta)$ is the zero fiber of the shifted
gradient map
$\widehat{\mu_{\liep^\beta}}=\mu_{\liep^\beta}-\beta$.

For $\beta\in\liep$ we set
$S^{\beta_+}:=\SS_{G^\beta}(\m_{\liep^\beta}(\beta))(X^{\beta_+}_{\norm\beta^2})$,
i.e., $S^{\beta_+}$ is the set of $G^\beta$-semistable points in
$X^{\beta_+}_{\norm\beta^2}$ with respect to the shifted gradient
map $\mu_{\liep^\beta}-\beta$. The set $S_\beta:=G\cdot
S^{\beta_+}$ is called the \emph{pre-stratum associated with
$\beta$}.

\begin{remark}
  $S^{\beta_+}$ is a locally closed submanifold of $X$ since it is an open subset of
  $X^{\beta_+}$.
\end{remark}

\begin{remark} As we already noted (see
  Remark~\ref{remark:independent}), the set
  $X_{\norm{\beta}^2}^{\beta_+}$ only depends upon the group
  $\Gamma=\exp(\RR\beta)$. However, in general, $S^{\beta_+}$ depends
  upon $G^\beta$.
\end{remark}

Let $\cp$ denote the set of critical points of $\eta_{\liep}\colon
X\to \RR$, $\eta_\liep(x)=\frac 12 \norm {\mup(x)}^2$ and set
$\bp:=\mup(\cp)$. Since $\eta_\liep$ is $K$-invariant the sets
$\cp$ and $\bp$ are $K$-stable.

We will now formulate our main general results.

\begin{theorem}[(Slice Theorem for Pre-Strata)]
  \label{theorem:strataslice} Let $\beta\in\bp$ and let $S_\beta$ be the
  pre-stratum associated with $\beta$.
  \begin{enumerate}
  \item The pre-stratum $S_\beta$ is a locally closed submanifold of $X$ and
$S^{\beta_+}$ is a $G^{\beta_+}$-stable locally closed submanifold of $X$.\label{item:theorem:strataslice1}
  \item The
    natural map $\twist G{G^{\beta_+}}S^{\beta_+}\to S_\beta$ is a
    diffeomorphism.\label{item:theorem:strataslice2}
  \item The natural map $\twist K{K^{\beta}}S^{\beta_+}\to S_\beta$
    is a diffeomorphism.\label{item:theorem:strataslice3}
  \end{enumerate}
\end{theorem}

The following is an analogue of the Quotient Theorem~\ref{theorem:quotient}  for
semistable points in the context of pre-strata.  Here
$\mathcal M_{\liep}(K\cdot\beta):=\mu_\liep\inv(K\cdot\beta)$
plays the role of $\mathcal M_\liep$ in
Theorem~\ref{theorem:quotient} and the case where $\beta=0$ is
just Theorem~\ref{theorem:quotient}.

\begin{theorem}[(Quotient Theorem for
  Pre-Strata)]\label{theorem:strataquotient} Let $\beta\in\bp$ and let
  $S_\beta$ be the pre-stratum associated with $\beta$. Then the
 topological Hilbert quotient $\quot{S_\beta}{G}$
  exists. Let  $\pi_{S_\beta}\colon S_\beta\to\quot{S_\beta}{G}$ denote the
  quotient map.
  \begin{enumerate}
  \item Every fiber of $\pi_{S_\beta}$ contains a unique closed $G$-orbit.
Any other orbit in
    the fiber has strictly larger dimension.\label{item:theorem:strataquotient1}
  \item The closure of every $G$-orbit in a $\pi_{S_\beta}$-fiber
    contains the closed $G$-orbit.\label{item:theorem:strataquotient2}
  \item Every fiber of $\pi_{S_\beta}$ intersects
    $\mp(K\cdot\beta)\cap S_\beta$ in
    a unique $K$-orbit which lies in the unique closed $G$-orbit.\label{item:theorem:strataquotient3}
  \item The inclusion $\mp(K\cdot\beta)\cap S_\beta\hookrightarrow S_\beta$ induces a homeomorphism
    $(\mp(K\cdot\beta)\cap S_\beta)/K\cong \quot {S_\beta}G$.\label{item:theorem:strataquotient4}
  \end{enumerate}
\end{theorem}

Note that we can apply Theorem~\ref{theorem:quotient} to the set
of semistable points
$S^{\beta_+}=\SS_{G^\beta}(\m_{\liep^\beta}(\beta))(X^{\beta_+}_{\norm\beta^2})$.
In the following discussion it will turn out that the quotients
$\quot{S^{\beta_+}}{G^\beta}$ and $\quot{S_\beta}G$ are isomorphic
(Proof of Theorem~\ref{theorem:strataquotient}) and that
$\m_{\liep^\beta}(\beta)\cap S^{\beta_+} \subset \mp(K\cdot
\beta)\cap S_\beta$ (Lemma~\ref{lemma:mp} and
Proposition~\ref{proposition:prestratum}). We may summarize all
this in the following commutative diagram:
\[
\xymatrix{ & \m_{\liep^\beta}(\beta)\cap S^{\beta_+}
  \ar@{^{(}->}[ld]\ar@{^{(}->}[rr]\ar[dd] \ar@{^{(}->}[rd]& &
  \mp(K\cdot \beta)\cap S_\beta \ar@{^{(}->}[ld]\ar[dd]
  \\
  S^{\beta_+}\ \ar@{^{(}->}[rr] \ar[dd] & & S_\beta \ar[dd] &
  \\
  & (\m_{\liep^\beta}(\beta)\cap S^{\beta_+})/K^\beta \ar@{=}[rr]^\sim
  \ar@{=}[rd]^\sim \ar@{=}[ld]_\sim& & (\mp(K\cdot\beta)\cap S_\beta)/K
  \ar@{=}[ld]_\sim
  \\
  \quot{S^{\beta_+}}{G^\beta}\ar@{=}[rr]^\sim & & \quot{S_\beta}{G} &
}
\]
We will also see that
$\m_{\liep^\beta}(\beta)\cap S^{\beta_+}$
is the set of minima of $\eta_\liep\vert S^{\beta_+}$ (Lemma~\ref{lemma:SbetaPlusMinimality}) and that
$\mp(K\cdot \beta)\cap S_\beta$ is the set of minima of $\eta_\liep\vert S_\beta$ (Proposition~\ref{proposition:minimality}).

\begin{theorem} [(Pre-Stratification Theorem)]\label{theorem:prestratification}
  Let $\beta$, $\tilde\beta\in\bp$.
  \begin{enumerate}
  \item If $S_{\tilde\beta}\cap \overline{S_\beta}\ne\emptyset$ and
    $K\cdot\beta\neq K\cdot\tilde\beta$, then
    $\|\tilde\beta\|>\|\beta\|$.\label{item:theorem:prestratification1}
  \item The following are equivalent.
    \begin{enumerate}
    \item $K\cdot \beta=K\cdot\tilde\beta$.\label{item:theorem:prestratification21}
    \item $S_{\beta}= S_{\tilde\beta}$.\label{item:theorem:prestratification22}
    \item $S_{\beta}\cap S_{\tilde\beta}\neq\emptyset$.\label{item:theorem:prestratification23}
    \end{enumerate}\label{item:theorem:prestratification2}
  \end{enumerate}
\end{theorem}

Now we will give additional information about the pre-strata which will lead to the proofs of the theorems.

\begin{lemma} \label{lemma:mp} For any $\beta\in\liep$ we have
  \begin{enumerate}
  \item $\mp(\beta)\cap X^\beta=\mp(\beta)\cap \cp$.\label{item:lemma:mp1}
  \item $\mp(\beta)\cap X^\beta=\m_{\liep^\beta}(\beta)\cap X^\beta$.\label{item:lemma:mp2}
  \item $\m_{\liep^\beta}(\beta)\cap X^\beta=\m_{\liep^\beta}(\beta)\cap
  X^{\beta_+}_{\norm\beta^2}$.\label{item:lemma:mp3}
  \end{enumerate}
\end{lemma}

\begin{proof}
The first part is a direct consequence of
Lemma~\ref{lemma:gradients}. Since $\mu$ is
$\exp(\im\er\beta)$-equivariant, on $X^\beta$ it takes values in
$(\lieu^*)^\beta$ and $\mup$ takes values in $\liep^\beta$. This
implies (\ref{item:lemma:mp2}). Part (\ref{item:lemma:mp3}) follows from
Corollary~\ref{corollary:hyperplane}.
\end{proof}

\begin{remark}
  Lemma~\ref{lemma:mp} implies that $S^{\beta_+}$ (equivalently $S_\beta$) is nonempty if and
  only if $\beta\in\bp$. This is the only reason why Theorem~\ref{theorem:strataslice} and Theorem~\ref{theorem:strataquotient} are
  formulated for $\beta\in\mathcal B_\liep$ and not for
  $\beta\in\liep$.
\end{remark}

As in Section~\ref{section:parabolic} we let 
$p^{\beta_+}\colon X_{\norm\beta^2}^{\beta_+}\to X_{\norm\beta^2}^\beta$ denote the
$G^\beta$-equivariant map $p^{\beta_+}(y)=\lim_{t\to-\infty}\exp
t\beta\cdot y$.

\begin{proposition}\label{proposition:prestratum}
  For every $\beta\in \bp$ we have the following.
  \begin{enumerate}
  \item
    $S^{\beta_+}=(p^{\beta_+})\inv(\SS_{G^\beta}(\m_{\liep^\beta}(\beta))(X^\beta))$.\label{item:proposition:prestratum1}
  \item $S^{\beta_+}$ is $G^{\beta_+}$-stable.\label{item:proposition:prestratum2}
  \item $\m_{\liep^\beta}(\beta)\cap
  S^{\beta_+}=\m_{\liep^\beta}(\beta)\cap X^\beta$.\label{item:proposition:prestratum3}
  \item
    $\SS_{G^\beta}(\m_{\liep^\beta}(\beta))(X_{\norm{\beta}^2}^\beta)$ is a
    $G^\beta$-equivariant strong deformation retract of $S^{\beta_+}$.\label{item:proposition:prestratum4}
  \end{enumerate}
\end{proposition}

\proof Let $y\in X^{\beta_+}$ and set $x=p^{\beta_+}(y)$. Since
$p^{\beta_+}$ is $G^\beta$-equivariant and fixes $X^\beta$ and
since $x$ lies in $\overline{G^\beta\cdot y}$,  we have that
$\emptyset\neq\overline{G^\beta\cdot x}\cap
X^\beta\cap\m_{\liep^\beta}(\beta)$ if and only if
$\emptyset\neq\overline{G^\beta\cdot y}\cap
X^\beta\cap\m_{\liep^\beta}(\beta)$.
 Hence
the equality in (\ref{item:proposition:prestratum1}) follows from Lemma~\ref{lemma:mp}. Invariance
of $S^{\beta_+}$ with respect to the $G^{\beta_+}$-action follows
from (\ref{item:proposition:prestratum1}) and (\ref{item:proposition:prestratum3}) follows from Lemma~\ref{lemma:mp} (\ref{item:lemma:mp3}). The deformation of
$S^{\beta_+}$ onto 
$\SS_{G^\beta}(\m_{\liep^\beta}(\beta))(X_{\norm{\beta}^2}^\beta)$ is given by 
$(t,y)\mapsto \exp t\beta\cdot y$.  \qed

\begin{corollary} \label{corollary:parabolicinvariance} For all
  $\beta\in \bp$ we have $S_\beta=K\cdot S^{\beta_+}$.
\end{corollary}

\begin{proof}
This follows from $G=KG^{\beta_+}$
(Lemma~\ref{lemma:GgleichKGbeta}), $G^{\beta_+}$-stability of
$S^{\beta_+}$ and the definition $S_\beta=G\cdot S^{\beta_+}$.
\end{proof}

By Theorem~\ref{theorem:quotient} the topological Hilbert
quotients $\quot{\mathcal S_{G^\beta}(\mathcal
M_{\liep^\beta}(\beta))(X^\beta)}{G^\beta}$ and
$\quot{S^{\beta_+}}{G^\beta}$ exist.
Proposition~\ref{proposition:prestratum} implies that the
inclusion $\mathcal S_{G^\beta}(\mathcal
M_{\liep^\beta}(\beta))(X^\beta)\hookrightarrow S^{\beta_+}$
induces an isomorphism of these quotients. Note that for $x,y\in
S^{\beta_+}$ we have $\overline{G^\beta\cdot
x}\cap\overline{G^\beta \cdot y}\cap S^{\beta_+} \neq\emptyset$ if
and only if $\overline{G^{\beta_+}\cdot
x}\cap\overline{G^{\beta_+}\cdot y}\cap S^{\beta_+} \neq\emptyset$
by Proposition~\ref{proposition:strata}. Therefore the quotient
$\quot{S^{\beta_+}}{G^{\beta_+}}$ exists and is not only
isomorphic to but also equal to $\quot{S^{\beta_+}}{G^\beta}$. We
summarize this discussion in the following commutative diagram.
\[
\xymatrix{
\mathcal S_{G^\beta}(\mathcal M_{\liep^\beta}(\beta))(X^\beta)\ \ar[d] \ar@{^{(}->}[r]    & S^{\beta_+}\ar[r]^{\mathrm{id}_{S^{\beta_+}}}\ar[d] &S^{\beta_+}\ar[d]\\
\quot{\mathcal S_{G^\beta}(\mathcal M_{\liep^\beta}(\beta))(X^\beta)}{G^\beta}\ar@{=}[r]^{\hspace{1.2cm}\sim} & \quot{S^{\beta_+}}{G^\beta}\ar@{=}[r]&\quot{S^{\beta_+}}{G^{\beta_+}}
}
\]

Recall that by Theorem~\ref{theorem:quotient} each fiber of the quotient $S^{\beta_+}\to\quot{S^{\beta_+}}{G^\beta}$ contains a unique closed $G^\beta$-orbit which is the unique orbit of minimal dimension in that fiber. We now show the analogous  fact for the action of $G^{\beta_+}$ on $S^{\beta_+}$.
We make use of the following remark.

\begin{remark} Let $H$ be a Lie group acting on a manifold $Y$. Then for all $d\in \N$, $\{y\in Y\mid \dim H\cdot y\geq d\}$ is open.

\end{remark}
\begin{proposition}\label{proposition:ParabQuotient}
Let $G^\beta\cdot x$ be the unique closed $G^\beta$-orbit in a fiber
$F$ of the quotient $S^{\beta_+}\to\quot{S^{\beta_+}}{G^\beta}$. Then
$G^{\beta_+}\cdot x$ is the unique closed $G^{\beta_+}$-orbit in $F$
and every other $G^{\beta_+}$-orbit in $F$ has strictly larger
dimension.
\end{proposition}

\begin{proof} Since $p^{\beta_+}(x)\in G^\beta\cdot x$, we must have
that $x\in X^\beta$. To simplify notation let $Q$ denote
$G^{\beta_+}$. Let $z\in F$ and set $y:=p^{\beta_+}(z)$. Then
$y\in F\cap X^\beta$.

Assume that $Q\cdot z\neq Q\cdot y$. We show that $\dim Q\cdot
z>\dim Q\cdot y$.
 Let $\lieq_y^c$ be a subspace of $\lieq$ complementary to
$\lieq_y$. Since an
$\exp(\er\beta)$-invariant neighborhood of $y$ in $X^{\beta_+}$
can be identified with an $\exp(\er\beta)$-invariant neighborhood
of $0$ in the tangent space $T_y(X^{\beta_+})$ (see
Section~\ref{section:parabolic}) and since $\lieq_y^c\cdot
y=\lieq\cdot y$ is $\exp(\er\beta)$-stable, there is
an $\exp(\er\beta)$-invariant locally closed submanifold $N_y$ of
$X^{\beta_+}$ with $y\in N_y$ and $T_y(X^{\beta_+})=\lieq_y^c\cdot y\oplus
T_y(N_y)$. This implies that the map $\Phi_y\colon\lieq_y^c\times
N_y\to X^{\beta_+}$,
$(\xi, w)\mapsto\exp(\xi)\cdot w$, is a local diffeomorphism near
$(0,y)$.  Since
$Q\cdot z$ intersects every open neighborhood of $y$ and since
$\exp(\lieq_y^c)\subset Q$, we may assume that there is a $z'\in Q\cdot z$ such that $z'\in N_y$ and such that
$\dim\lieq_y^c\cdot z'=\dim\lieq_y^c\cdot y$. Since $Q\cdot
z\neq Q\cdot y$, $z'\not\in X^\beta$ and $0\neq\beta_X(z')\in T_{z'}(N_y)$ does not lie in $\lieq_y^c\cdot z'$. Thus $\dim Q\cdot z'\geq\dim\lieq^c_y+1>\dim Q\cdot y$.

It now suffices to show that $\dim Q\cdot y>\dim Q\cdot x$ if
$Q\cdot y\neq Q\cdot x$.  We have that  $\dim G^\beta\cdot x<\dim
G^\beta\cdot y$. Since $p^{\beta_+}(x)=x$   the
$G^\beta$-equivariance and $R^{\beta_+}$-invariance of
$p^{\beta_+}$ imply that $\dim Q\cdot x=\dim G^\beta\cdot x+\dim
R^{\beta_+}\cdot x$ and similarly for $y$. Since we may assume
that $y$ is arbitrarily close to $x$, we have that $\dim
R^{\beta_+}\cdot y\geq  \dim R^{\beta_+}\cdot x$. Thus $\dim
Q\cdot y= \dim G^\beta\cdot y+\dim R^{\beta_+}\cdot y>\dim
G^\beta\cdot x+\dim R^{\beta_+}\cdot x=\dim Q\cdot x$.

We have shown that the orbit $Q\cdot x$ is of minimum dimension in
$F$. It remains to show that it is the unique closed orbit. Suppose that $x'\in\overline{Q\cdot x}$ where $Q\cdot x'\neq Q\cdot x$. Then $\dim Q\cdot x\geq\dim Q\cdot x'$. But this contradicts the fact that
$\dim Q\cdot x<\dim Q\cdot x'$. Hence $Q\cdot x$ is closed.
By the definition of topological Hilbert
quotient there are no other closed orbits in the fiber $F$.
\end{proof}

\begin{lemma}\label{lemma:SbetaPlusMinimality}
Let $\beta\in\bp$ and $z\in \overline{S^{\beta_+}}$. Then $\eta_{\liep}(z)\ge\frac{1}{2}\|\beta\|^2$ and equality holds if and only if $z\in\mp(\beta)\cap X^\beta$.
\end{lemma}

\begin{proof}
We have $\norm{\mup(z)}\ge\norm{\mu_{\liep^\beta}(z)}\ge\norm\beta$
  by Corollary~\ref{corollary:hyperplane}. If equality holds then $\mu_{\liep^\beta}(z)=\beta$ and $z$ is a minimum of $\eta_{\liep^\beta}\vert G^\beta z$. Therefore it is a critical point of $\eta_{\liep^\beta}$ and we have $z\in\m_{\liep^\beta}(\beta)\cap X^\beta=\m_{\liep}(\beta)\cap X^\beta$ by Corollary~\ref{corollary:critical}.
\end{proof}

\begin{proposition} \label{proposition:minimality}
Let $\beta\in\bp$ and let $Y_\beta$ be a $G$-stable closed subset of $S_\beta$ and
    $\overline{Y_\beta}$ its closure in $X$.
 We have
    $\eta_\liep\vert\overline{Y_\beta}\ge\frac 12\norm\beta^2$ and the
    set of $z\in \overline{Y_\beta}$ where $\eta_\liep(z)=\frac
    12\norm\beta^2$ coincides with $\mp(K\cdot \beta)\cap Y_\beta$
    and is non-empty if $Y_\beta$ is non-empty.
\end{proposition}

\begin{proof}
We have $\overline{Y_\beta}=K\cdot \overline{Y^{\beta_+}}$ where
  $Y^{\beta_+}:=Y_\beta\cap S^{\beta_+}$
  (Corollary~\ref{corollary:parabolicinvariance}). For
  $z\in\overline{Y_\beta}$ we have $z=k\cdot y$ where $k\in K$ and $y\in
  \overline{Y^{\beta_+}}\subset \overline{S^{\beta_+}}$. Then
  $\norm{\mup(z)}=\norm{\mup(y)}\ge\norm\beta$. If we have
 equality then $\mu_{\liep^\beta}(y)=\beta$ and $y\in X^\beta$ (Lemma~\ref{lemma:SbetaPlusMinimality}).
 Now $Y^{\beta_+}$ is closed in $S^{\beta_+}$. Therefore $y\in \m_\liep(\beta)\cap Y^{\beta_+}$ and $z\in \mp(K\cdot\beta)\cap Y_\beta$.

 Finally, $\mp(K\cdot \beta)\cap Y_\beta$ is non-empty since
 $\overline{G^\beta \cdot y}$ intersects $\m_{\liep^\beta}(\beta)\cap
X^\beta$ for every $y\in Y^{\beta_+}\subset S^{\beta_+}$.
\end{proof}

\begin{corollary}\label{corollary:GBahnAbschlussMinimal}
The closure of every $G$-orbit in $S_\beta$ intersects
$\m_\liep(\beta)\cap X^\beta$.
\end{corollary}

We now work towards the proof of the Slice Theorem for Pre-Strata.

\begin{lemma}\label{lemma:parabolicslice}
  Let $\beta\in\bp$. The identity $S^{\beta_+}\to S^{\beta_+}$ induces a diffeomorphism
  $\phi\colon\twist K{K^\beta}{S^{\beta_+}}\to \twist
  G{G^{\beta_+}}{S^{\beta_+}}$.
\end{lemma}

\begin{proof}
 Since $G=KG^{\beta_+}$, $K\cap G^{\beta_+}=K^\beta$ and
  $\lieg=\liek+\lieg^{\beta_+}$, $\phi$ is a 1-1 onto submersion between manifolds of the same dimension, hence a diffeomorphism.
\end{proof}

In the following we fix $\beta\in\bp$ and define $\Psi\colon \twist
G{G^{\beta_+}}{S^{\beta_+}}\to S_\beta$ by $\Psi([g,s])=g\cdot s$.

\begin{lemma} \label{lemma:immersion} The map $\Psi$ is an immersion.
\end{lemma}
\begin{proof}
  First we show that $\Psi$ is an immersion at any $x\in\mp(\beta)\cap
  S^{\beta_+}$.
 Since $\lieg=\liek+\lieg^{\beta_+}$ and
  $\lieg^{\beta_+}\cdot x \subset T_x(S^{\beta_+})$ it is sufficient
  to show that $\zeta\in\liek$ and $\zeta_X(x)\in T_x(S^{\beta_+})$
  implies that $\zeta\in \liek^\beta\subset \lieg^{\beta_+}$.  This is
  a consequence of Lemma~\ref{lemma:hessian2} since $T_x(S^{\beta_+})$
  is the sum of the eigenspaces with non-negative eigenvalues of the
  isotropy action of $\beta$ on $T_x(X)$.
  Since every $G^\beta$-orbit in ${S^{\beta_+}}$ intersects every open
  neighborhood of $\mathcal{M}_{\liep}(\beta)\cap S^{\beta_+} $
  and $\Psi$ is $G$-equivariant it follows that $\Psi$ is an
  immersion.
\end{proof}

\begin{proposition}\label{proposition:limit}
  Let $w_n\in \twist G{G^{\beta_+}}{S^{\beta_+}}$ be a sequence such
  that
  $\Psi(w_n)$ converges to $x\in\mp(\beta)\cap
  X^\beta$. Then the sequence $w_n$ has a convergent subsequence and
  every convergent subsequence converges to $[e,x]$.
\end{proposition}

\begin{proof}
  Let $w_n=[k_n,s_n]$ where $k_n\in K$ and $s_n\in S^{\beta_+}$
  (Lemma~\ref{lemma:parabolicslice}). We may assume that $k_n$
  converges to $k\in K$. This implies that $\lim_{n\to
    \infty}\Psi(k_n\inv[k_n,s_n])=\lim_{n\to \infty} s_n=k\inv\cdot
  x\in\overline{S^{\beta_+}} \subset
  \overline{X^{\beta_+}_{\norm\beta^2}}$. From $\norm{\mup(k\inv\cdot
    x)}=\norm{\mup(x)}=\norm\beta$ and $k\inv\cdot x\in
  \overline{S^{\beta_+}}$ we obtain that $k\inv x\in \mp(\beta)\cap X^\beta$ (Lemma~\ref{lemma:SbetaPlusMinimality}).
 Since also $x\in\mp(\beta)\cap X^\beta$ we conclude that  $k\in K^\beta$ and that $\lim_{n\to\infty}[k_n,s_n]=[k,k\inv\cdot x]=[e,x]$.
\end{proof}

\begin{corollary} \label{corollary:injective} The map $\Psi$ is
  injective.
\end{corollary}
\begin{proof}
  We first show that for every $x\in\mp(\beta)\cap X^\beta$ there
  exists an open neighborhood $\Omega$ of $x\in X$ such that
  $\Psi\vert \Psi\inv(S_\beta\cap \Omega)$ is injective.  If
  $\Psi\vert\Psi\inv(S_\beta\cap\Omega)$ is not injective for $\Omega$
  sufficiently small, then there are $v_n,w_n\in \twist
  G{G^{\beta_+}}{S^{\beta_+}}$ such that $v_n\ne w_n$, $\Psi(v_n)=\Psi(w_n)$ and $\lim_{n\to
    \infty} \Psi(v_n)=x$. Since $\Psi$ is
  an immersion (Lemma~\ref{lemma:immersion}) its restriction to some
  open neighborhood $W$ of $[e,x]$ in $\twist
  G{G^{\beta_+}}{S^{\beta_+}}$ is injective. We may assume that
  $\lim_{n\to \infty} v_n=\lim_{n\to \infty} w_n=[e,x]$
  (Proposition~\ref{proposition:limit}).  This contradicts injectivity
  of $\Psi\vert W$.

Since every $G$-orbit in $S_\beta$ intersects each neighborhood of a point of  $\mp(\beta)\cap X^\beta$  (Corollary~\ref{corollary:GBahnAbschlussMinimal}) and since $\Psi$ is $G$-equivariant,  every $G$-orbit in $S_\beta$  intersects an open subset $\Omega$ of $X$ non trivially such that
  $\Psi\vert\Psi\inv(S_\beta\cap \Omega)$ is injective. Hence $\Psi$ is injective.
 \end{proof}

\begin{proof}[\textbf{Proof of the Slice Theorem for Pre-Strata
    \ref{theorem:strataslice}}] We first show (\ref{item:theorem:strataslice1}). It follows from Proposition \ref{proposition:strata} that $S^{\beta_+}$ is a locally closed submanifold of $X$. Let $x\in \mathcal
  M_\liep(\beta)\cap X^\beta$ and choose an open neighborhood $W_x$ of
  $[e,x]\in \twist G{G^{\beta_+}}{S^{\beta_+}}$ such that $\Psi(W)$ is
  a submanifold of $X$ (Lemma~\ref{lemma:immersion}). Then there is an
  open neighborhood $\Omega_x$ of $x$ in $X$ such that $\Omega_x\cap
  S_\beta\subset \Psi(W)$ since otherwise we could use
  Proposition~\ref{proposition:limit} to arrive at a contradiction.
  The open sets $g\cdot \Omega_x$ cover $S_\beta$
  (Corollary~\ref{corollary:GBahnAbschlussMinimal}) and the open sets $g\cdot W_x$ cover $ \twist
  G{G^{\beta_+}}{S^{\beta_+}}$ where $g\in G$, $x\in \mathcal
  M_\liep(\beta)\cap X^\beta$. Since $g\cdot \Omega_x\cap
  S_\beta\subset \Psi(g\cdot W_x)=g\cdot \Psi(W_x)$ we have (\ref{item:theorem:strataslice1}).

  Part (\ref{item:theorem:strataslice2}) follows from (\ref{item:theorem:strataslice1}), Lemma~\ref{lemma:immersion} and
  Corollary~\ref{corollary:injective}.
  Lemma~\ref{lemma:parabolicslice} and (\ref{item:theorem:strataslice2}) imply (\ref{item:theorem:strataslice3}).
\end{proof}

\begin{proof}[\textbf{Proof of the Quotient Theorem for Pre-Strata
    \ref{theorem:strataquotient}}]
Recall that the topological Hilbert quotients
$\quot{S^{\beta_+}}{G^\beta}$ and
$\quot{S^{\beta_+}}{G^{\beta_+}}$ exist and coincide. The Slice
Theorem for Pre-Strata implies that $\quot{S_\beta}G$ exists and
that the inclusion $S^{\beta_+}\subset S_\beta$ induces an
isomorphism $\quot{S^{\beta_+}}{G^{\beta_+}}\cong\quot{S_\beta}G$.

By Proposition~\ref{proposition:ParabQuotient} every fiber of the quotient map
$\pi_{S^{\beta_+}}\colon S^{\beta_+}\to \quot{S^{\beta_+}}
{G^{\beta_+}}$ contains a unique closed $G^{\beta_+}$-orbit $G^{\beta_+}\cdot
x$ and every other orbit in that fiber has strictly larger
dimension. Assertion (\ref{item:theorem:strataquotient1}) then
follows from the Slice Theorem for Pre-Strata.

Assertion (\ref{item:theorem:strataquotient2}) follows from the definition of the topological
Hilbert quotient and (\ref{item:theorem:strataquotient1}).

  Let $q\in \quot{S_\beta}{G}$ and  set $F_q:=(\pi_{S_\beta})\inv(q)$. Then
  $F_q\cap S^{\beta_+}$ is a fiber of $\pi_{S^{\beta_+}}$ and intersects
  $\m_{\liep^\beta}(\beta)$ in a unique $K^\beta$-orbit $K^\beta\cdot x$ where
  $x\in \m_{\liep^\beta}(\beta)\cap X^\beta$ and $G^\beta\cdot x$ is the unique
  closed orbit in $F_q\cap S^{\beta_+}$. Hence $G\cdot x$ is the unique closed orbit
  in $F_q$ and $K\cdot x$ is the intersection of $F_q$ with $\mp(K\cdot\beta)$ and we have (\ref{item:theorem:strataquotient3}).

 The assertion in  (\ref{item:theorem:strataquotient4}) follows from (\ref{item:theorem:strataquotient1}),  (\ref{item:theorem:strataquotient3}) and Theorem  \ref{theorem:quotient}.
\end{proof}

\begin{proof}[\textbf{Proof of the Pre-Stratification Theorem
    \ref{theorem:prestratification}}] Let $y\in
  \overline{S_{\beta}}\cap S_{\tilde\beta}$. Then $\overline{G \cdot
    y}$ contains a point $x\in\mp(\tilde\beta)$
  (Corollary~\ref{corollary:GBahnAbschlussMinimal} applied to $\tilde \beta$)
  and is contained in $\overline{S_{\beta}}$. This implies that
  $\|\tilde\beta\|\ge\|\beta\|$ and equality holds if and only if
  $K\cdot \tilde\beta=K\cdot \beta$
  (Proposition~\ref{proposition:minimality}).
  Hence we must have $\|\tilde\beta\|>\|\beta\|$, proving (\ref{item:theorem:prestratification1}).

  We now show (\ref{item:theorem:prestratification2}). Assume that (\ref{item:theorem:prestratification21}) holds and choose $k\in K$ such
that $\tilde\beta=k\cdot \beta$. Then $S^{\tilde\beta_+}=k\cdot S^{\beta_+}$
implies that $S_{\tilde\beta}=S_\beta$. Hence (\ref{item:theorem:prestratification21}) implies (\ref{item:theorem:prestratification22}).
Obviously (\ref{item:theorem:prestratification22}) implies (\ref{item:theorem:prestratification23}). Assume (\ref{item:theorem:prestratification23}) and let $z\in S_\beta\cap S_{\tilde\beta}$. If (\ref{item:theorem:prestratification21}) fails, then (\ref{item:theorem:prestratification1}) shows that $\|\tilde\beta\|>\|\beta\|$ and that $\|\beta\|>\|\tilde\beta\|$, a contradiction.
\end{proof}

\section{Consequences and special cases}
\label{section:applications}

In this section we point out special cases and several
consequences of the results obtained in the previous section.

First note that for $\beta=0$ the pre-stratum $S_\beta$ coincides
with $S^{\beta_+}$ and is the set of semistable points
$\SS_G(\m_\liep)$. The Pre-Strata Slice Theorem
\ref{theorem:strataslice} is trivial in this case and the Quotient
Theorem~\ref{theorem:strataquotient} is just
Theorem~\ref{theorem:quotient}. If $G=G^\beta$, in particular if
$G$ is commutative, then $S_\beta=S^{\beta_+}$ coincides with the
set of semistable points
$\SS_G(\m_{\liep}(\beta))(X^{\beta_+}_{\|\beta\|^2})$.

For Morse theoretical considerations the following should be noted.

\begin{lemma} \label{lemma:stratacriticalpoints} Let $\beta\in
  \bp$. Then \( \cp\cap S_\beta=\mp(K\cdot\beta)\cap S_\beta.  \)
\end{lemma}

\begin{proof}
 We have $\mp(K\cdot\beta)\cap S_\beta= K\cdot (\mp(\beta)\cap X^\beta) \subset\cp$ by Corollary \ref{corollary:critical}.  Conversely, if
  $x\in \cp\cap S_\beta$,
   let $\tilde\beta:=\mu_\liep(x)$. Then $x\in \mp(\tilde\beta)\cap X^{\tilde\beta}\subset S_{\tilde\beta}$. Since $x\in S_\beta\cap S_{\tilde\beta}$, we have $K\cdot\beta=K\cdot\tilde\beta$, and  hence
   $x\in \mp(K\cdot
  \beta)$.
\end{proof}

\begin{corollary}
  If $X=\SS_G(\mp)$, then $x\in X$ is a critical point of
  $\eta_{\liep}$ if and only if $\eta_{\liep}(x)=0$, i.e., the
  critical set of $\eta_{\liep}$ coincides with $\mp$.
\end{corollary}

Fix $\beta\in\bp$.
The first two theorems in Section~\ref{section:strata}
contain precise information  about the structure of the $G$-action
on $S_\beta$ and the quotient map $\pi_{S_\beta}\colon
S_\beta\to \quot{S_\beta}{G}$.
We now wish to apply the Slice Theorem for the action of $G^\beta$ on $S^{\beta_+}$. To simplify notation, let us assume that $X=S_\beta$ and let $Y$ denote $S^{\beta_+}$.
We have an artificial $G^{\beta_+}$-action on $Y^\beta$ defined by
$(g,y)\mapsto g\acts y:=\pi^{\beta_+}(g)\cdot y$, $y\in Y^\beta$, $g\in G^{\beta_+}$.
The mapping $p^{\beta_+}\colon Y\to Y^\beta$ is
equivariant with respect to the given action on $Y$ and the artificial
action on $Y^\beta$. Here we use the notation
introduced in Section~\ref{section:parabolic}.
Note that the $G^{\beta_+}$-isotropy group at $y\in Y^\beta$ is   the semi-direct
product $H_y:=(G^\beta)_y\ltimes R^{\beta_+}$. Applying the Slice Theorem we obtain the following.
\begin{proposition}
Let $x\in\m_{\liep^\beta}(\beta)\cap X^\beta\subset Y^\beta$. Then there is a $(G^\beta)_x$-stable submanifold $Y_0$ of $Y^\beta$
containing $x$ such that
$G^\beta\cdot Y_0=G^{\beta_+}\acts Y_0$ is open in $Y^\beta$ and such that the
natural map
$\twist{G^{\beta_+}}{H_x}{Y_0}\to G^{\beta_+}\acts Y_0$
is an isomorphism. Let $Y_1:=(p^{\beta_+})\inv(Y_0)$. Then $Y_1$ is
an $H_x$-stable submanifold of $Y$,
\[
\twist{G^{\beta_+}}{H_x}{Y_1}\to G^{\beta_+}\cdot Y_1
\]
is an isomorphism
and $G^{\beta_+}\cdot Y_1$ is open in $Y$. Moreover,  $G\cdot Y_1$ is
open in $X$ and
\[
\twist{G}{G^{\beta_+}}({\twist{G^{\beta_+}}{H_x}{Y_1}})\to G\cdot Y_1
\]
is a $G$-equivariant isomorphism.
\end{proposition}
 \begin{remark}
The entire construction can be carried out in a manner compatible with the
various Hilbert quotients associated to the $G$-action on $X$. That is, we can arrange that  $G^\beta\cdot Y_0$ is saturated with respect to the quotient map
$Y^\beta\to \quot{Y^\beta}{G^\beta}$ and that $G^\beta\cdot Y_1$ is saturated with respect
to $Y\to \quot{Y}{G^\beta}$, etc.
\end{remark}
\begin{remark}
The fiber $F$  of $\pi_{G^\beta}\colon Y\to \quot{Y}{G^\beta}$ through $x\in \m_{\liep^\beta}(\beta)\cap Y$
is $G^\beta$-isomorphic to $\twist{G^\beta}{(G^\beta)_x}  F_x$
where $F_x$ is a closed $(G^\beta)_x$-stable subset of
the $(G^\beta)_x$-stable subspace $W^{\beta_+}$ (Section~\ref{section:parabolic})
of $T_x(X)$.
In this identification
every $(G^\beta)_x$-orbit in $F_x$ contains $0\in W^{\beta_+}$ in its closure.
Since $F$ is also a fiber of $\pi_{G^{\beta_+}}$, the $(G^\beta)_x$-space $F_x$
is equipped with a  topological action of $H_x=(G^\beta)_x\ltimes R^{\beta_+}$.
\end{remark}

Before we discuss open pre-strata for general $X$ we need to compute the
Hessian of $\eta_\liep$ at critical points.

\begin{proposition}\label{proposition:indexhessian} Let $x\in\cp$  be
   a critical point of $\eta_{\liep}\colon X\to \RR$
 and let $S_\beta$
  be the associated pre-stratum. Let $H_x(\eta_{\liep})$ denote the
  Hessian of $\eta_{\liep}$ at $x$. Then
  \begin{enumerate}
  \item $H_x(\eta_{\liep})= 0$ on $T_x(K\cdot x)$,\label{item:kx}
  \item $H_x(\eta_{\liep})> 0$ on $\liep^\beta\cdot
  x+\lier^{\beta_+}\cdot x$, where $\lier^{\beta_+}$ is the Lie algebra of $R^{\beta_+}$,\label{item:prx}
  \item $H_x(\eta_{\liep})\ge 0$ on $T_x(S_\beta)=\lieg\cdot
    x+T_x(S^{\beta_+}) =\liek\cdot x+T_x(S^{\beta_+})$ and
  \item $H_x(\eta_{\liep})< 0$ on $T_x(S_\beta)^\perp=(\lieg\cdot
    x)^\perp\cap T_x(S^{\beta_+})^\perp= (\liek\cdot x)^\perp\cap
    T_x(S^{\beta_+})^\perp$.
  \end{enumerate}
\end{proposition}

\begin{proof} $K$-invariance of $\eta_\liep$ implies (\ref{item:kx}). The other assertions follow from
$ (\lieg\cdot x)^\perp\subset (\liep\cdot x)^\perp=\ker d\mup(x)$,
Proposition~\ref{proposition:hessian}, Proposition~\ref{proposition:strata} (\ref{item:proposition:strata3}) and the fact that
$d\beta_X(x)\colon T_xX\to T_xX$ has strictly positive eigenvalues on $\lier^{\beta_+}\cdot x$.
\end{proof}

\begin{remark}
We have a tangent space decomposition
$T_x(G\cdot x)=T_x(K\cdot x)\oplus\liep^\beta\cdot x\oplus\lier^{\beta_+}\cdot x$. This follows from the
decompositions $G=K\cdot G^{\beta_+}$, $G^{\beta_+}=G^\beta\cdot
R^{\beta_+}$, the identity $K\cap G^{\beta_+}=K^\beta$ and the fact that $G^\beta$ acts on $X^\beta$ whereas
$R^{\beta_+}$ acts on the fibers of $p^{\beta_+}$ (see Section~\ref{section:parabolic}).
Thus the behavior of $H_x(\eta_\liep)$ on $T_x(G\cdot x)$ is precisely
described by Proposition~\ref{proposition:indexhessian},
(\ref{item:kx}) and (\ref{item:prx}).
\end{remark}

\begin{corollary}\label{corollary:minimalvalue}
  If $\beta\in\bp$ is such that $\lVert\beta\rVert^2$ is a minimum
  value of $2\eta_{\liep}$, then the corresponding pre-stratum
  $S_\beta$ is open in $X$ and coincides with
  $\mathcal{S}_G(\mp(\beta))$.
\end{corollary}
\proof Let $x\in\cp$ such that $\mup(x)=\beta$. By
Proposition~\ref{proposition:indexhessian}, $S_\beta$ is open in $X$.
It contains $\m_{\liep^\beta}(\beta)\cap X^\beta=\mp(\beta)\cap
X^\beta=\mp(\beta)$ since $\norm\beta^2$ is a minimum. But
$\mathcal{S}_G(\mp(\beta))$ is the smallest $G$-stable open set
containing $\mp(\beta)$, so that $\mathcal{S}_G(\mp(\beta))\subset
S_\beta$. On the other hand,
\[
S_\beta=G\cdot S_{G^\beta}(\m_{\liep^\beta}(\beta)\cap
X^\beta)(X^{\beta_+})\subset S_G(\mp(\beta)).
\]
\qed

The argument in the proof shows  the following more general result.

\begin{corollary}\label{corolarry:locallyminimalvalue}
  Let $x$ be a local minimum of $\eta_\liep$. Then $x\in\cp$ and the union
  of the connected components of the corresponding stratum $S_\beta$
  which intersect $K\cdot x$ non-trivially is an open subset of $X$.
\end{corollary}

\begin{corollary}\label{corollary:GHomogen} Suppose that $X=G\cdot x$ where $x\in\cp$ and
  $\beta:=\mup(x)$. Then
  \begin{enumerate}
    \item $G^{\beta_+}\cdot x=S^{\beta_+}$.\label{item:corollary:GHomogen1}
      \item $\SS_{G^\beta}(\m_{\liep^\beta}( \beta))(X^\beta)=G^\beta\cdot
    x$.\label{item:corollary:GHomogen2}
  \item $K^\beta\cdot x=\mp(\beta)$.\label{item:corollary:GHomogen3}
  \item $\cp=K\cdot x$.\label{item:corollary:GHomogen4}
  \end{enumerate}
\end{corollary}

\proof In the homogeneous case the Slice Theorem for Pre-Strata
gives us a $G$-equivariant isomorphism
$\twist{G}{G^{\beta_+}}{S^{\beta_+}}\to G\cdot x$ and therefore a
$G$-equivariant map $p\colon G\cdot x\to G/G^{\beta_+}$ with
$G^{\beta_+}\cdot x=p\inv(G^{\beta_+})=S^{\beta_+}$. This shows
(\ref{item:corollary:GHomogen1}).

Let $z=q\cdot x\in S^{\beta_+}$ where $q\in G^{\beta_+}$. Then
$\lim_{t\to -\infty}\exp(t\beta) \cdot q\cdot x=h\cdot x$ where
$h=\lim_{t\to -\infty} \exp (t\beta) q \exp(-t\beta)\in G^\beta$. Hence
we have  (\ref{item:corollary:GHomogen2}).

>From (\ref{item:corollary:GHomogen2}) we obtain that
$\m_{\liep^\beta}(\beta)\cap X^\beta=K^\beta\cdot x$. This implies
(\ref{item:corollary:GHomogen3}). Now
(\ref{item:corollary:GHomogen3}) gives $\mp(K\cdot \beta)=K\cdot
x$ and (\ref{item:corollary:GHomogen4}) follows from
Lemma~\ref{lemma:stratacriticalpoints}.  \qed

\begin{corollary}
  Let $x\in\cp$. Then $x$ is a global minimum of $\eta_\liep|_{G\cdot x}$.
\end{corollary}

\begin{corollary}\label{corollary:homogeneouslocalmaximum}
  If $\eta_{\liep}\vert G\cdot x$ has a local maximum at $x$, then
  $G\cdot x=K\cdot x$.
\end{corollary}

\begin{proof} If $x$ is a local maximum of $\eta_{\liep}$, then it is a
critical point of $\eta_{\liep}$. But every critical point of
$\eta_{\liep}$ on $G\cdot x$ is a global minimum.  This implies that
the set of critical points, which is a $K$-orbit, is open in $G\cdot
x$. Since $K$ intersects every component of $G$, we have $K\cdot
x=G\cdot x$.
\end{proof}

\begin{corollary}\label{corollary:Gorbit=Korbit}
If a $G$-orbit $G\cdot x$ in $X$ is compact then $G\cdot x=K\cdot x$.
\end{corollary}

In Corollary~\ref{corollary:homogeneouslocalmaximum} we have a
special case of the existence of a maximal pre-stratum $S_\beta$.
By this we mean that $\norm\beta\ge\norm{\mup(z)}$ holds for every
$z\in X$. In this case $\eta_\liep\vert S_\beta$ is constant and
$\mup$ maps $S_\beta$ equivariantly onto $K\cdot \beta$. This
implies that every $G$-orbit in $S_\beta$ is a $K$-orbit, hence
closed. More precisely, we have $S_\beta=\mp(K\cdot \beta)\cap
S_\beta$ and $S^{\beta_+}=\m_{\liep^\beta}(\beta)\cap
S^{\beta_+}=\mp(\beta)\cap X^\beta$. Note that this implies that
$S^{\beta_+}$ coincides with a union of connected components of
$X^\beta$. Also since $S_\beta\cong \twist
{K}{K^\beta}{(\m_{\liep^\beta}(\beta)\cap X^\beta)}$ we have that
$S_\beta$ is compact if $X^\beta$ is compact or $\mup$ is a proper
map.\\
\\
As an application, we show how our results imply Matsuki duality.
\begin{corollary}\label{corollary:matsuki}
Assume that $X=Z$ is $U$-homogeneous, i.e., a generalized flag
manifold, and assume that $G$ is a real form of $U^\CC$. Then each
$G$-orbit and each $K^\ce$-orbit intersect  $\cp$ in a unique
$K$-orbit (Matsuki duality). If $x\in\cp$, then $K\cdot x$ is the
set of global minima of $\eta_\liep$ on $G\cdot x$ and of
$\eta_\liek$ on $K^\ce\cdot x$. Thus $G\cdot x\cap K^\ce\cdot
x=K\cdot x$. Moreover $K^\ce\cdot x$ (resp. $G\cdot x$) is open in
$X$ if and only if $G\cdot x=K\cdot x$ (resp. $K^\ce\cdot x=K\cdot
x$).
\end{corollary}

\begin{proof}
We have $T_xX=T_x(G\cdot x)+T_x(K^\ce\cdot x)$. We may assume that
the decomposition  $\lieu^\CC=\liek\oplus \im\liek\oplus\im\liep
\oplus \liep$ is orthogonal with respect to a $U$-invariant inner
product on $\lieu^\CC$. Then
$\frac{1}{2}\|\mu\|^2=\eta_\liek+\eta_\liep$ where
$\eta_\liek:=\frac{1}{2}\|\mu_\liek\|^2$. Since $Z$ is
$U$-homogeneous this function is constant and $\eta_\liek$ is the
negative of $\eta_\liep$ up to a constant. In particular, the set
of critical points of $\eta_\liek$ equals $\cp$. Using the
behaviour of the Hessian of $\eta_\liep$ (resp. $\eta_\liek$)
described in Proposition~\ref{proposition:indexhessian}, we get
that each $G$-orbit (resp. $K^\ce$-orbit) is open in the
$G$-pre-stratum (resp. $K^\ce$-pre-stratum) in which it is
contained. Therefore each $G$-orbit and each $K^\ce$-orbit
intersects $\cp$ in a unique $K$-orbit.  if $x\in\cp$ then $K\cdot
x$ is the set of global minima of $\eta_\liep$ on $G\cdot x$ and
of $\eta_\liek$ on $K^\ce\cdot x$. Again, since
$\frac{1}{2}\|\mu\|^2$ is constant, we have $G\cdot x\cap
K^\ce\cdot x=K\cdot x$. The last statement follows from the
identities $G\cdot x\cap K^\ce\cdot x=K\cdot x$ and
$T_xX=T_x(G\cdot x)+T_x(K^\ce\cdot x)$.
\end{proof}

\begin{remark} Matsuki duality does not depend upon the choice of the moment map, as follows. Assume that we are given another  $U$-equivariant moment map $\tilde\mu\colon Z\to\lieu^*$.
Let $\tilde\eta_\liep:=\frac{1}{2}\|\tilde\mu_\liep\|^2$ and let $\tilde{\mathcal C}_\liep$
denote the set of critical points of $\tilde\eta_\liep$. If there is an $x_0\in\cp$ where
$x_0\notin\tilde{\mathcal C}_\liep$, then, by Matsuki duality, there is an
$x_1\in K^\ce\cdot x_0\cap\tilde{\mathcal C}_\liep$ with $x_1\notin\cp$ and
$\eta_\liep(x_0)>\eta_\liep(x_{1})$. Again, by Matsuki duality, we obtain an $x_2\in G\cdot x_1\cap\cp$
with $x_2\notin\tilde{\mathcal C}_\liep$ and $\eta_\liep(x_1)>\eta_\liep(x_2)$.
Inductively we obtain a sequence $x_n$ in $\cp$ with $\eta_\liep(x_n)>\eta_\liep(x_{n+1})$.
But $\cp$ consists of finitely many $K$-orbits so $K\cdot x_{n_0}=K\cdot x_0$ for
some $n_0>0$, which contradicts the inequality $\eta_\liep(x_0)>\eta_\liep(x_{n_0})$. Hence $\tilde{\mathcal C}_\liep=\cp$.
\end{remark}

We can also easily obtain an old result of Wolf \cite{Wolf}.

\begin{corollary} \label{corollary:old wolf}Assume that $X=Z$ is $U$-homogeneous
and  that $G$ is a real form of $U^\CC$. Then there is exactly one
closed $G$-orbit in $X$.
\end{corollary}

\begin{proof}
If $G\cdot x$ is a closed orbit in $Z$ it equals $K\cdot x$ by
Corollary~\ref{corollary:Gorbit=Korbit}. Then we may assume
$x\in\cp$. By Corollary~\ref{corollary:matsuki} the dual orbit
$K^\ce \cdot x$ is open. Since $K^\CC$ is complex and acts
holomorphically on $Z$ there is only one open orbit.
\end{proof}

\section{Proper moment maps}\label{section:proper}

In this section we consider what happens when $\mup\colon X\to
\liep$ is a proper map. In this case $\eta_\liep\colon X\to \RR$
is a proper non-negative function, i.e., an exhaustion of $X$.
Note that properness of $\eta_\liep$ implies that $X$ is a closed
$G$-stable submanifold of $Z$.

Let $(t,x)\mapsto \varphi_t(x)$ denote the local flow of the vector field
$\grad\eta_{\liep}$ on $X$.  Properness of $\mup$ implies that
$\varphi_t$ exists for all negative $t$.  By
Lemma~\ref{lemma:gradients}, the flow $\varphi_t$ is along $G$-orbits.
Note that $\eta_{\liep}(\varphi_t(x))$ is strictly increasing unless
$x$ is a critical point of $\eta_{\liep}$. For $z\in X$ let $\mathcal
L(z)$ denote the set of points $y\in X$ such that $y=\lim_{k\to\infty}
\varphi_{t_k}(z)$ for some sequence $t_k $ which goes to $-\infty$.
Since we assume $\mup$ to be proper we have the following.

\begin{lemma}\label{lemma:gradientlimit}
  Let $z\in X$. Then $\mathcal L(z)$ is non empty and for every open
  neighborhood $V$ of $\mathcal L(z)$ there exists a $t_0\in\RR$ such
  that $\varphi_t(z)\in V$ for every $t<t_0$.
\end{lemma}

\begin{proof}
  Let $I_V:=\{t\in\RR: t\le 0\text{ and $\varphi_t(z)\in X\setminus
    V$}\}$. Since $V$ contains $\mathcal L(z)$ the properness of
  $\eta_\liep$ implies that $I_V$ is compact. For $t_0:=\min\{t\in
  I_V\}$ we have $\varphi_t(z)\in V$ for all $t<t_0$.
\end{proof}

For $\beta\in\bp$ let $\mathcal C_\liep(K\cdot\beta)$ denote
$\cp\cap\mp(K\cdot\beta)$ and set
$\cp(\norm\beta)=\cp\cap\eta_\liep\inv(\frac 12\norm\beta^2)$.  If
we set $\beta:=\mup(x)$ for some $x\in \mathcal L(z)$, then
$x\in\cp(K\cdot \beta)=\cp(K\cdot \beta)\cap S_\beta\subset
\mp(K\cdot \beta)\cap S_\beta$. It is a consequence of
Lemma~\ref{lemma:gradientlimit} that $\mathcal L(z)$ is connected.
Then $\eta_\liep$ is constant on
  $\mathcal{L}(z)$ and we have
  \[
  \mathcal{L}(z)\subset
  \bigcup_{\norm{\tilde\beta}=\norm\beta}\cp(K\cdot
  \tilde\beta)=\cp(\norm\beta).
  \]

We need the following technical result whose proof we give at the
end of this section.

\begin{proposition} \label{proposition:normaldistance} Let $z\in X$
  and choose $\beta\in\bp$ such that $\mathcal{L}(z)\cap\cp(K\cdot
  \beta)\not=\emptyset$. Then there is an open neighborhood $\Omega$
  of $ \cp(K\cdot \beta)$ and a smooth function
  $\rho\colon \Omega\to \RR$ such that
  \begin{enumerate}
  \item $\rho\ge 0$,
  \item $\Omega\cap S_\beta=\{z\in \Omega: \rho(z)=0\}$ and
  \item $d\rho(\grad(\eta_\liep))(z)\le 0$ for all $z\in\Omega$.
  \end{enumerate}
\end{proposition}

The proposition has the following consequence for the flow
$\varphi_t$.

\begin{theorem}[(Stratification Theorem)]\label{theorem:stratification}
  If $\mup\colon X\to \liep$ is proper, then $X=\cup_\beta S_\beta$
  where $\beta$ runs through a complete set of representatives of
  $K$-orbits in $\bp$. The union is disjoint, each $S_\beta$ is a
  locally closed submanifold of $X$ and
  \[\overline{S_\beta}\subset S_\beta\cup
  \{S_\gamma\colon\norm\gamma>\norm\beta\}.\]
\end{theorem}

\begin{proof}
Let $z\in X$ and $\beta\in\bp$ such that $\mathcal{L}(z)\cap
  \cp(K\cdot\beta)\neq \emptyset$. Choose an open neighborhood $\Omega$
  of $\cp(K\cdot\beta)$ which has the properties given in
  Proposition~\ref{proposition:normaldistance}.

   Let $\Omega_0$ be an open neighborhood of $\cp(K\cdot \beta)$ which
  is relatively compact in $\Omega$ and let $Q:=\Omega_0\cap
  \cp(\norm\beta)$. Since $\cp(\norm\beta)\setminus Q$ is compact and
  $(\cp(\norm\beta)\setminus Q)\cap \overline{S_\beta}=\emptyset$
  there is an open relatively compact neighborhood $V'$ of
  $\cp(\norm\beta)\setminus Q$ in $X$ such that $\overline{V'}\cap
  \overline{S_\beta}=\emptyset$. Then the set $V:=V'\cup \Omega_0$ is
  an open neighborhood of $\mathcal{L}(z)\subset\cp(\norm\beta)$.  Let $r:=\min\{\rho(y):
  y\in \overline{V'}\cap \overline{\Omega_0}\}$. Since
  $\overline{V'}\cap S_\beta=\emptyset$ we have $\rho(y)>0$ on the
  compact set $\overline{V'}\cap\overline{\Omega_0}$ and therefore
  $r>0$.  The set $\Omega_1:=\{y\in\Omega_0: \rho(y)<r\}$ is an open
  neighborhood of $\cp(K\cdot \beta)$ such that $\overline{\Omega_1}\cap
  V\subset \Omega$.

  Since $\mathcal{L}(z)\cap \cp(K\cdot\beta)\neq\emptyset$ there exists a $t_0$ with $\varphi_{t_0}(z)\in\Omega_1$. By Lemma~\ref{lemma:gradientlimit} we may assume that $\varphi_t(z)\in V$ for all $t\le t_0$. We claim that $\varphi_t(z)\in\Omega_1$ for all $t\le t_0$. It suffices to show that $I_{\Omega_1}:=\{t\in\RR: t\le t_0,\ \varphi_t(z)\in\Omega_1\}$ is
  connected  since it is not bounded from below.
If this would not be the case then we would find a
  connected component $I_1=(a_1,t_1)=\{t\in I_{\Omega_1}: a_1<t<t_1\}$ of
  $I_{\Omega_1}$ where $a_1$ is possibly $-\infty$ and $t_1<t_0$.  Then
  $\varphi_{t_1}(z)\in V\cap (\overline{\Omega_1}\setminus \Omega_1)$, i.e.,
  $\rho(\varphi_{t_1}(z))=r$, and $\rho(\varphi_t(z))< r$ for all
  $t\in I_1$. This contradicts $a_1<t_1<t_0$, since $t\mapsto
  \rho(\varphi_{-t}(z))$ is defined and increasing for all $t$ in a
  sufficiently small open interval which contains $t_1$.

  Now assume that $z\not\in S_\beta$. Since
  $\rho(\varphi_t(z))>0$ for all $t<t_0$ and $t\mapsto
  \rho(\varphi_{-t}(z))$ is increasing there is no $x\in
  \mathcal{L}(z)\cap\cp(K\cdot\beta)\subset S_\beta$ with
  $\rho(x)=0$. This contradiction proves that $z\in S_\beta$ and shows that
  $X=\bigcup_\beta S_\beta$. The union is disjoint and we also have
  $\overline{S_\beta}\subset S_\beta\cup \{S_\gamma\colon
  \norm\gamma>\norm\beta\}$ by the Pre-Stratification
  Theorem~\ref{theorem:prestratification}.
\end{proof}

In the proof we have seen that $z\in S_\beta$ if $\mathcal
L(z)\cap\cp(K\cdot \beta)\neq\emptyset$. Since $\mathcal L(z)\subset\cp(\|\beta\|)$ we obtain

\begin{corollary}
Let $z\in S_\beta$. Then $\mathcal L(z)\subset\cp(K\cdot\beta)$.
\end{corollary}

\begin{remark}
  The map $\varphi\colon (-\infty,0]\times S_\beta\to S_\beta$, $(t,x)\mapsto
  \varphi_t(x)$ extends to a continuous map
  $\bar\varphi\colon
  [-\infty,0]\times S_\beta\to S_\beta$. Then $\bar\varphi$ is a
  $K$-equivariant strong deformation retraction of $S_\beta$ onto
  $\mathcal C_\liep(K\cdot\beta)$ which stabilizes
  closures of $G$-orbits. In particular, the set $\mathcal L(z)$
  consists of one point. For the proof one uses a result of Marle
  (\cite{Marle}) and Guillemin and Sternberg (\cite{GS}) which states
  that there exist local coordinates at any $z\in Z$ in which
  the moment map $\mu$ is real analytic. Applying this to points
  $x\in\mathcal C_\liep(K\cdot\beta)$
  one obtains the existence of $\bar\varphi$ as in
  \cite{Ne85} (see \cite{Sch89})
  using an inequality of \L ojasiewicz.

If $\mu_\liep$ is not assumed to be proper, one can prove that
$\bar\varphi$ realizes
  $\mathcal M_\liep$ as a strong deformation retract of a neighborhood
  of $\mathcal M_\liep$ in $S_0=\mathcal S_G(\mathcal M_\liep)$. Here
  the same method applies since in \cite{HSt05} it is shown that
there exist relatively compact neighborhoods of points in
$\mathcal M_\liep$ which have the property that they contain
$\varphi_t(z)$ for all $t<t_0$ if they contain $\varphi_{t_0}(z)$.
\end{remark}

We have the following characterization of $S_\beta$.

\begin{corollary}\label{corollary:Strataproper}
  Let $\mup$ be proper and $\beta\in\bp$. Then
  \begin{align*}
    S_\beta&=\{z\in X\colon \beta\in \mup(\overline{G\cdot z}) \text{
      and } \lVert\beta\rVert\le\lVert\mu_\liep(g\cdot z)\rVert \text{
      for all } g\in
    G\}\\
    &=\{z\in X\colon
    \lim_{t\to-\infty}\varphi_t(z)\in\cp(K\cdot\beta)\}.
  \end{align*}
\end{corollary}

\begin{remark}
  In general we have $S_\beta\subset\{z\in X\colon \beta\in
  \mup(\overline{G\cdot z}) \text{ and }
  \lVert\beta\rVert\le\lVert\mu_\liep(g\cdot z)\rVert \text{ for all }
  g\in G\}$.
\end{remark}

We now work towards a proof of
Proposition~\ref{proposition:normaldistance}.

As in Section 2, if $x$ is a zero of a vector field $\zeta$ on
$X$, then $d\zeta(x)$ denotes the corresponding endomorphism of
$T_x(X)$.

\begin{lemma}\label{lemma:gradientlinearization}
  For $x\in \cp$ and $\beta:=\mup(x)$ we have
  \begin{enumerate}
  \item $d\grad\eta_\liep(x)\cdot v=(d\beta_X(x))\cdot v +
    (d\mup(x)\cdot v)_X(x)$ for all $v\in T_x(X)$ and \label{item:lemma:gradientlinearization1}\item the linear
    map $L\colon T_x(X)\to T_x(X)$, $v\mapsto (d\mu_\liep(x)\cdot
    v)_X(x)$ maps $T_x(S_\beta)$ into itself and is zero on
    $T_x(S_\beta)^\perp$.\label{item:lemma:gradientlinearization2}
  \end{enumerate}
\end{lemma}

\begin{proof}
  Let $\gamma$ be a smooth curve through $x\in \cp$ such that
  $\gamma(0)=x$ and set $v:=\frac{d}{dt}\vert_0\gamma(t)$. We have
  $\grad \eta_\liep (\gamma(t))= \mup(\gamma(t))_X(\gamma(t))=\grad
  \mup^{\mup(\gamma(t))}(\gamma(t))$ and $\mup(\gamma(t))=\mup(x)+ t\cdot
  d\mup(x)\cdot v +t^2 R(t)$ for a continuous map $R(t)$. Since
  $\grad\eta_\liep(\gamma(t))=\grad\mup^{\mup(x)}(\gamma(t))+t
  \grad\mup^{d \mup(x)\cdot
    v}(\gamma(t))+t^2\grad\mu_\liep^{R(t)}(\gamma(t))$ we have
  $d\grad\eta_\liep(x)\cdot v=d\beta_X(x)\cdot v+(d\mup(x)\cdot
    v)_X(x)$ proving (\ref{item:lemma:gradientlinearization1}). We also have (\ref{item:lemma:gradientlinearization2}), since $\lieg\cdot
  x\subset T_x(S_\beta)$ and $\ker d\mup(x)=(\liep\cdot x)^\perp$.
\end{proof}

\begin{proof}[Proof of Proposition~\ref{proposition:normaldistance} ]
Let $\Omega_0$ be an open neighborhood of $S_\beta$ in $X$ which
is diffeomorphic to the normal bundle $T(S_\beta)^\perp\vert
S_\beta$ of $S_\beta$ in $X$ such that $S_\beta$ corresponds to
the zero section. Using this identification we define
$\rho\colon\Omega_0\to \er$, $\rho(v)=(v,v)_{\pi(v)}$ where
$\pi\colon T(S_\beta)^\perp\vert S_\beta\to S_\beta$ denotes the
canonical projection and $(\cdot,\cdot)_{\pi(v)}$ denotes the
inner product on $T_{\pi(v)}(X)$. Then we have $\rho\ge 0$ and
$\Omega_0\cap S_\beta=\{y\in\Omega_0: \rho(y)=0\}$.

For $x_0\in\cp(K\cdot\beta)$ there exists an open neighborhood
$U_0$ of $x_0$ in $\Omega_0$ which can be identified with
$\er^{n}\times\er^m$ where $x_0$ corresponds to $0$, $S_\beta$
corresponds to $\er^n\times\{0\}$ and the normal space
$T_x(S_\beta)^\perp$ corresponds to $\{x\}\times\er^m$ for $x\in U_0\cap
S_\beta$. Furthermore, this identification can be chosen such that
$\rho$ has the form $\rho(v_1,v_2)=\|v_2\|^2$ on
$\er^n\times\er^m$.

Expressing the content of Lemma~\ref{lemma:gradientlinearization}
in these coordinates we obtain
$$
\grad \eta_\liep(v_1,v_2)=d
\beta_X(0)(v_1,v_2)+L(v_1)+\text{O}(\lVert v_1\rVert^2+\lVert
v_2\rVert^2)
$$
where $L$ is a linear mapping from $\RR^{n}$ to
$\RR^n$.

Since $d \beta_X(0)$ has only negative eigenvalues on $\RR^m$, the
quadratic form $v_2\mapsto (v_2, d {\beta}_X(0)\cdot v_2)$ is
strictly negative definite. Furthermore, $d\beta_X(0)$ preserves
$T_0(S_\beta)$ and
  $T_0(S_\beta)^\perp$, so we have
$$
d\rho (\grad\eta_\liep)(v_1,v_2)=(v_2, d\beta_X(0)\cdot v_2)+
\text{O}(\lVert v_1\rVert^3+\lVert v_2\rVert^3).
$$
Thus after possibly shrinking $U_0$ we can find a $c<0$ such that
$d\rho(\grad\eta_{\liep})(v_1,v_2)\leq c\cdot\rho(v_1,v_2)$ for
$(v_1,v_2)\in U_0$.

Since $x_0\in\cp(K\cdot\beta)$ was arbitrary we find a neighborhood
$\Omega\subset\Omega_0$ of $\cp(K\cdot\beta)$ in $X$ such that
$d\rho(\grad\eta_\liep)(z)\le 0$ if
$z\in\Omega$.
\end{proof}

\section{Compact manifolds}\label{section:compact}

In this section we show that for a compact $G$-stable submanifold
$X$ of $Z$ there are only finitely many pre-strata. Let $\liea$ be
a maximal subalgebra of $\liep$ and $A=\exp\liea$ the
corresponding subgroup of $G$. Note that $A$ is compatible with
the Cartan decomposition of $U\c$. Since the corresponding
$A$-gradient map $\mua$ is locally constant on the smooth compact
manifold $X^A$ we find that $\mathcal V:=\mua(X^A)$ consists of
finitely many points.

In \cite{HSt05} the following is established.

\begin{lemma} \label{lemma:convexity} For $x\in X$ the image
  $\mu_{\mathfrak{a}}(A\cdot x)$ is an open convex subset of the
  affine subspace $\mu_{\mathfrak{a}}(x)+\mathfrak{a}_x^\perp$ of
  $\mathfrak{a}$ where $\mathfrak{a}_x^\perp:=\{\beta\in\liea;\
  \langle\alpha,\beta\rangle=0 \text{ for all } \alpha\in\liea_x\}$ denotes the orthogonal complement of
  $\mathfrak{a}_x$ in $\mathfrak{a}$.
\end{lemma}

\begin{proposition}\label{prop:KonvexeHuelleVonFixpunkten}
  Let $C\subset X$ be compact and $A$-invariant. If the image
  $\mu_\liea(C)$ is convex, then it is the convex hull of a subset of
  $\mathcal V$. In particular, it is a polytope.
\end{proposition}

\begin{proof}
  Since $\mu_\liea(C)$ is compact and convex it is the convex hull of
  the set of its extreme points.
 Let $\beta$ be an extreme point and let
  $x\in C$ with $\mu_\liea(x)=\beta$. Since $C$ contains $A\cdot x$,
  Lemma~\ref{lemma:convexity} implies that $\beta$ lies in an open
  subset of $\beta+\liea_x^\perp$ which is contained in
  $\mu_\liea(C)$. But $\beta$ is an extreme point so
  $\liea_x^\perp=\{0\}$. Therefore $\liea_x=\liea$ and consequently
  $x\in X^A$.
\end{proof}

\begin{corollary}
  Let $X$ be compact and $x\in X$. Then the image
  $\mu_\liea(\overline{A\cdot x})$ is the convex hull of a subset of
  $\mathcal V$.
\end{corollary}

\begin{corollary}\label{corollary:bpFinite} Let $X$ be compact and
  $\beta\in\bp\cap\liea$.  Then $\beta$ is the closest point to zero
  in the convex hull of finitely many elements of $\mathcal V$. In
  particular, $\bp\cap\liea$ is a finite set and $\bp$ consists of
  finitely many $K$-orbits.
\end{corollary}

\begin{proof}
  Let $x\in\cp$ such that $\mup(x)=\beta$. Note that $x$ is also a
  critical point for $\eta_{\liea}$. Applying Proposition~\ref{proposition:minimality} to the group $A$ and $Y_\beta:=\overline{A\cdot x}$,
  we obtain that $\norm\beta\le \norm \zeta$ for all $\zeta\in\mua(Y_\beta)$.
  Since $\mua(Y_\beta)$ is the convex hull of images of $A$-fixed points,
  $\beta$ is the closest point to the origin of a subset of $\mathcal
  V$.
\end{proof}

\begin{corollary} If $X$ is compact, then there are only finitely
  many pre-strata.
\end{corollary}

In the case where $X$ is compact we call a pre-stratum $S_\beta$
a stratum.

\section{Morse inequalities}\label{section:morse}

Let $G$ be an arbitrary Lie group and $Y$ a topological $G$-space.
The $G$-equivariant cohomology $\eco G Y$ of $Y$ is by definition
the ordinary cohomology $\co{\operatorname{EG}\times^GY}$ where
$\operatorname{EG}$ is the total space of a universal $G$-bundle
$p\colon \operatorname{EG}\to \operatorname{BG}$ and
$\operatorname{EG}\times^GY$ is the quotient of
$\operatorname{EG}\times Y$ by the diagonal action of $G$.

We define the $G$-equivariant Poincar\'{e} series to be the power
series $\poin GY:=\sum_{n=0}^\infty t^n\dim\eco GY$ where
cohomology is computed with respect to a field $\mathbb{K}$ which
we omit in the notation. If $G=\{e\}$ is the trivial group,
equivariant cohomology coincides with ordinary cohomology and we
write $\poin{}{Y}$ for $\poin GY$.

Now assume that $X$ is an orientable
 connected compact $G$-invariant submanifold
of $Z$ and that $G$ is a closed connected compatible subgroup of
$U^\ce$. Consider the finite decomposition $X=\cup_\beta S_\beta$
(Theorem~\ref{theorem:stratification}, Corollary~\ref{corollary:bpFinite}). If $S_{\beta,m}$ denotes
the union of the connected components of $S_\beta$ which have codimension
$m$ in $X$, then $X=\cup_{\beta, m}S_{\beta,m}$ is again a decomposition of $X$ with   closures of strata having the properties analogous to those in Theorem \ref{theorem:stratification}. In the following, cohomology
is computed with coefficients in $\qu$ if all the strata $S_\beta$ are
orientable submanifolds of $X$ and with coefficients in
$\zet_2=\zet/2$ otherwise.

Standard considerations in equivariant topology show that we have the following.

\begin{morseineq}There exists a series $R(t)$
  with non-negative integer coefficients such that
  \[\sum_{\beta, m}t^{m}\poin K {S_{\beta, m}}-\poin
  KX=(1+t)R(t)\] where $\beta$ runs through a complete set of
  representatives of $K$-orbits in $\bp$.
\end{morseineq}

Since $K$ is a strong deformation retract of $G$ and the strata are $G$-stable the Morse inequalities are also valid if one replaces $K$-equivariant co\-ho\-mo\-lo\-gy by $G$-equivariant cohomology.

The following result of \cite{Kirwan} is useful for computing equivariant cohomology.
\\
\\
\textit{Let $Z$ be compact and $K$ compact and connected. Then $\poin
  KZ=\poin{} Z\cdot\poin{}{\operatorname{BK}}$ where
  $\operatorname{BK}$ is the base space of a universal $K$-bundle.}
\\
\\
In \cite{Kirwan} it is shown that in the case where $X:=Z$ is
compact and $G=K^\ce$ is complex reductive the Morse inequalities
are equalities, i.e., $R(t)=0$. Furthermore, the strata $S_\beta$
are complex and therefore orientable in this case. In our general
situation where $G$ is not complex reductive, there may exist
non-orientable strata.

\begin{example}
  Let $X:=Z:=\pzweice$ be equipped with the standard action of
  $\Sldreice$ and with the moment
  map induced by the Fubini-Study
  metric. Then for $G=\Sldreier$ there are two strata, namely
  $\pzweier$, which is not orientable, and its complement.
\end{example}

We end with an example which shows that, in our situation,
the
Morse inequalities are not necessarily equalities, even in the case that
$X=Z$.

\begin{example}
  Let $U^\ce=\slzweice$ act on $X:=Z:=\peinsce$ by $g\cdot [z]:=[gz]$.
  Then the moment
  map with respect to the Fubini-Study metric is
  given by $\mu^\xi([z])=\frac{<\xi\cdot z,z>}{\im<z,z>}$ where
  $<\cdot,\cdot>$ is the standard Hermitian product on $\ce^2$. Let
  $G=\slzweier$. Then $K=\sozweier\cong S^1$ and $\liep$ is the set of
  real symmetric matrices of trace zero. One can show that $\peinsce$
  decomposes into two orientable strata, namely
  $\SS_G(\m_\liep)=G\cdot[1,\im]\cup G\cdot [1,-\im]$ and $S:=G\cdot
  [1,0]=K\cdot[1,0]\cong S^1$.

  For $S^1$, a universal bundle is given by $p\colon
  S^\infty\to\mathbb{P}_\infty(\ce)$ where
  $S^\infty:=\lim_{n\to\infty}S^{2n+1}$,
  $\mathbb{P}_\infty(\ce):=\lim_{n\to\infty}\mathbb{P}_n(\ce)$ and $p$
  is induced by the projections $S^{2n+1}\to
  S^{2n+1}/S^1\cong\mathbb{P}_n(\ce)$ where $S^1$ acts by
  multiplication on $S^{2n+1}\subset\ce^{n+1}$.

  We have \[\poin
  K{\peinsce}=\poin{}{\peinsce}\cdot\poin{}{\operatorname{BK}}=(1+t^2)\sum_{n\geq0}t^{2n}=1+2\sum_{n\geq
    1}t^{2n}.\] Since $[1,\im]$ (resp. $[1,-\im]$) can be realized
  $K$-equivariantly as a strong deformation retract of $G\cdot[1,\im]$
  (resp. $G\cdot[1,-\im]$), we have
  \begin{align*} \poin K{\SS_G(\m_\liep)}&=\poin K
    {G\cdot[1,\im]}+\poin K{G\cdot
      [1,-\im]}=2\cdot\poin{}{EK\times_K{[1,\im]}}\\&=2\cdot\poin{}{BK}=2\cdot\sum_{n\ge
      0}t^{2n}.\end{align*} For the second stratum we get
  \[\poin K S=\poin{}{\operatorname{EK}}=1.\] Finally we get the
  Morse inequalities
  \begin{align*}
    t^0\cdot\poin K{\sgmip}+t^1\cdot\poin KS-\poin
    K{\peinsce}&=2\sum_{n\ge 0} t^{2n}+t-(1+2\sum_{n\geq
      1}t^{2n})\\
    &=1+t.
  \end{align*}
  Thus $R(t)=1\neq 0$.
\end{example}

\newcommand{\noopsort}[1]{} \newcommand{\printfirst}[2]{#1}
\newcommand{\singleletter}[1]{#1} \newcommand{\switchargs}[2]{#2#1}
\providecommand{\bysame}{\leavevmode\hbox to3em{\hrulefill}\thinspace}

\end{document}